\documentclass[journal]{IEEEtran}
\usepackage{graphicx}
\usepackage{multirow}
\usepackage{subfig}
\usepackage{algorithm}
\usepackage{algpseudocode}
\usepackage{amsmath}
\usepackage{float}
\usepackage{amssymb}
\usepackage{lipsum}
\usepackage{cite}
\usepackage{comment}
\makeatletter
\newcommand*{\rom}[1]{\expandafter\@slowromancap\romannumeral #1@}
\makeatother

\begin{document}
	
\title{A Stochastic Programming Approach for Electric Vehicle Charging Network Design}

\author{Sina Faridimehr, Saravanan Venkatachalam, Ratna Babu Chinnam
\thanks{S. Faridimehr, S. Venkatachalam, and R. Chinnam are with the Department of Industrial and Systems Engineering at Wayne State University, Detroit, Michigan (e-mail: sina.faridimehr@wayne.edu; saravanan.v@wayne.edu;
	ratna.chinnam@wayne.edu), Corresponding author: S.Venkatachalam}}

\maketitle

\begin{abstract}
	Advantages of electric vehicles (EV) include reduction of greenhouse gas and other emissions, energy security, and fuel economy. The societal benefits of large-scale adoption of EVs cannot be realized without adequate deployment of publicly accessible charging stations. We propose a two-stage stochastic programming model to determine the optimal network of charging stations for a community considering uncertainties in arrival and dwell time of vehicles, battery state of charge of arriving vehicles, walkable range and charging preferences of drivers, demand during weekdays and weekends, and rate of adoption of EVs within a community. We conducted studies using sample average approximation (SAA) method which asymptotically converges to an optimal solution for a two-stage stochastic problem, however it is computationally expensive for large-scale instances. Therefore, we developed a heuristic to produce near to optimal solutions quickly for our data instances. We conducted computational experiments using various publicly available data sources, and benefits of the solutions are evaluated both quantitatively and qualitatively for a given community.
\end{abstract}

\begin{IEEEkeywords}
	two-stage stochastic programming, electric vehicle, charging network, sample average approximation.
\end{IEEEkeywords}

\IEEEpeerreviewmaketitle

\section{Introduction}

Electric vehicles (EVs) hold much promise including diversification of the transportation energy feedstock, reduction of greenhouse gas and other emissions, and improved public health by improving local air quality. In general, widespread adoption of EVs is in alignment with sustainable transportation objectives due to its social, economic, and environmental perspectives. It is estimated that an EV that draws its power from the U.S. electrical grid emits at least 30\% less $CO_2$ than comparable gasoline or diesel-fueled vehicles \cite{berger2015comparison}. As EV usage for daily commute increases, the consideration for the ability to recharge these vehicles away from home will become even more important. Ever-growing need to recharge EVs away from home necessitates designing effective networks of charging stations. Using multiple linear regression, Sierzchula et al. \cite{sierzchula2014influence} examined the effect of consumer financial incentives and several socio-economic factors on national EV market shares of 30 countries for the year 2012. The analysis shows that installing one charging station (per 100,000 residents) could have twice the impact on EV adoption rate compared to a \$1,000 financial incentive.

Many studies have been done on locating charging stations for EVs. However, majority of them concentrated on large-scale state-wide networks and only a few articles have investigated design of public charging station network in an urban area. Existing papers on charging station location problem often assume that demands for charging service are deterministic and known to the decision makers, while in reality, the traffic flows are stochastic in nature (varying by hour of day, weekday, weekend, commute purpose, destination etc) and carry significant uncertainty. The optimal solution of a deterministic model might become infeasible and/or significantly sub-optimal in the presence of these uncertainties. This paper adds to the growing field of designing EV charging station network by proposing a two-stage stochastic programming model to determine location and size of charging stations for a community. Considering uncertainties in charging pattern, demand, and drivers' behavior, the proposed stochastic model provides more robust charging network design decisions and thus access to charging service can be improved. However, a two-stage stochastic programming model often needs a large number of scenarios for good representation of uncertainties. We use sample average approximation (SAA) method as this will asymptotically converge to an optimal solution for a two-stage stochastic problem. SAA is a Monte Carlo simulation-based sampling technique in which we approximate the expected value of the objective function using a finite sample of scenarios. Since SAA can only solve small size problems within reasonable amount of time in general, an effective heuristic is also proposed for large-scale instances. The two-stage model and solution approach are evaluated by a case study constructed using the data representing Detroit midtown area in Michigan, U.S. In summary, the major contributions of this paper include: (1) formulation of a two-stage stochastic programming model to determine the location and capacity of public EV charging stations in an urban area to maximize access; (2) incorporation of uncertainties in EV demand flows, EV drivers' charging patterns, arrival and departure time, purpose of arrival to a community, and preferred walking distance; (3) adoption of SAA to solve the two-stage model; (4) an effective heuristic that provides near optimal solutions for large-scale instances; and (5) a case study representing public charging network planning in Detroit midtown area and a post-analysis framework to analyze the outputs of the two-stage model on accessibility and utilization of charging service. The remainder of this paper is organized as follows: A review of related literature is presented in Section \ref{s2}. Section \ref{prob} provides problem description and the uncertainties considered in our model. Model formulation and the solution methodology are presented in Section \ref{model}. Section \ref{case} presents the case study, scenario construction, computational experiments and evaluations of results. Finally, conclusion and directions for future studies are provided in Section \ref{con}.

\section{Literature review}\label{s2}

During the last decade, many researchers have focused on optimally locating alternative-fuel-vehicle refueling stations. However, most of them are focused on EV charging network in large networks to cover demand between cities and metropolitan areas, and only a few articles examined the design of charging network in a community or an urban area. We review the existing literature related to design of an EV charging network and categorize it into two major groups: (A) deterministic approach which assumes that all parameters and demand are known for charging station network problem, and (B) stochastic approach that considers uncertainties regarding available budget for constructing charging network, type of charging stations, total short-term and long-term charging demand, and charging behavior of EV drivers. 

\subsection{Deterministic approach}

Upchurch et al. \cite{upchurch2009model} introduced capacitated flow refueling location model that considers a limit on the traffic flow that any location can refuel to maximize vehicle miles traveled by alternative-fuel vehicles. Frade et al. \cite{frade2011optimization} proposed a maximal covering model to find the optimal location of EV charging stations in an urban area by maximizing covered demand within a given distance. To deal with the computational burden of generating combinations of locations capable of serving the round trip on each route, a mixed-binary-integer optimization model is developed \cite{capar2012efficient}. Capar et al. \cite{capar2013arc} presented a more computationally efficient model for flow-refueling location model to answer some strategic questions such as what is the minimum number of charging stations required for refueling a certain percentage of traffic flow; and what are the impacts of refueling demand forecast on the location of fuel stations. A mixed-integer programming method to model capacitated multiple-recharging-station-location problem considering budget constraint and vehicle routing behavior, and using the concepts of set coverage and maximum coverage is proposed \cite{wang2013locating}. The model in \cite{baouche2014efficient} finds the optimal locations of charging stations for EVs in an urban area while minimizing total costs, consisting the travel cost from demand zones to charging locations and investment cost. Cavadas et al. \cite{cavadas2015mip} proposed a mixed-integer programming model to locate slow-charging stations for EVs in an urban environment considering the possibility that there might be several stops by each driver during the day and the driver can only charge the vehicle at one of these locations. Since tour-based network equilibrium model can precisely track the state-of-charge (SOC) of the battery and also consider the dwell time at each destination, model is proposed to optimally locate public charging stations for EVs considering recharging behavior of drivers \cite{he2015deploying}. Huang and Zhou \cite{huang2015optimization} developed an integer programming formulation to minimize the lifetime cost of equipment, installations, and operations of charging stations for plug-in EVs at workplaces by considering different charging levels and demographics of employees. In order to maximize the amount of vehicle-miles-traveled for an EV, a model is presented to select the optimal locations for public charging stations considering vehicle travel patterns \cite{shahraki2015optimal}. The authors applied their model on vehicle trajectory data of taxi fleet over a three week period in Beijing, China. A major limitation with all these studies is that they assume a deterministic problem setting. As we confirm through our experiments, employing a stochastic formulation can lead to a significant improvement in the objective of the planners.  

\subsection{Stochastic approach}

While planning under uncertainty has been addressed in many settings such as transportation, energy, disaster planning, supply chain management and production planning, the literature considering uncertainty in planning for EV charging network is limited. By considering both the transportation system and the power grid, Pan et al. \cite{pan2010locating} developed a two-stage stochastic programming model to find the optimal locations for battery exchange stations for plug-in hybrid electric vehicles (PHEV) accounting for uncertainty in demand for battery, loads, and generation capacity of renewable power sources. Tan and Lin \cite{tan2014stochastic} formulated the EV charging problem as a flow capturing location-allocation problem. They compared a deterministic case where charging demand is fixed over time to a stochastic one where consumer demand for charging service is random, and concluded that stochastic programming provides more realistic results. Hosseini and MirHassani \cite{hosseini2015refueling} proposed a two-stage stochastic program to locate permanent and portable charging stations with and without considering capacities to maximize the served traffic flows. A stochastic flow-capturing location model is also developed to locate a predetermined number of fast EV-charging stations within a given region considering uncertainties in EV flows \cite{wustochastic}.

To efficiently assist city planners and policy makers in planning for public EV charging network within a community, we need to adequately capture uncertainties that exists in demand for public charging service. To the best of our knowledge, this is the first study to address the problem of locating public EV charging stations for a community using a two-stage stochastic programming approach while accounting for uncertainties in total customer demand for public charging service, arrival and dwell time, battery SOC at the time of arrival, preference for charging away from home and willingness to walk patterns of EV drivers.

\begin{comment}
 In addition, we employ sample average approximation method to solve the stochastic programming problem and propose a heuristic to reduce the running time of the problem. We also evaluate the impact of public EV charging network planning on utilization level of charging stations and walking patterns within community through a case study. 
\end{comment}

\section{Problem Description and Uncertainties}\label{prob}

Unlike a conventional vehicle, an EV must often be parked for several hours to be recharged. Hence, public parking facilities are considered as potential locations for installing charging stations, which can in turn improve access to EVs as well as their adoption. Maximum number of installable charging stations depends on the total capacity of a parking lot. Without loss of generality, we assume that all charging point terminal types are semi-rapid charging ones (level 2 type charging stations) that are typically recommended for public and private parking lots, and provide 10 to 20 miles range per hour of charging. Also, a driver's walking distance to final destination is considered as the decisive contributing factor in choosing a parking lot \cite{bergantino2014drives}. Based on a driver's walking distance preference, we determine a possible set of parking lots that a driver can park the EV and then driver is randomly assigned to one of them. If charging stations are installed in any of parking lots that are within a driver's walking distance preference, driver will be attracted to one of those parking lots depending on the availability of a charging station at the time of arrival. If there is no parking lot within the maximum distance that a driver is willing to walk, we assume that driver will park the car on street, and since it is difficult to track the walking distance to final destination in this case, this demand is not considered in our analysis. It is also assumed that once a driver starts using a charging station, vehicle would not be unplugged until driver's activity is finished. 

Designing a public EV charging network entails estimation of demand for charging service. Like facility location models, we assume that demand occurs at fixed points on a network. Demand will be attracted to different parking lots based on drivers' willingness to walk to use charging stations. Scenarios representing demand uncertainty in the two-stage model will represent time and purpose of arrival to the community, EV's battery SOC at the time of arrival, duration of activity, drivers' preference for charging away from home and willingness to walk based on demographics, community size and seasonality factors. The following uncertainties are considered to affect demand for public EV charging stations:

\subsection{State of charge}

\begin{figure}[!htbp]
	\centering
	\includegraphics[scale = 0.3]{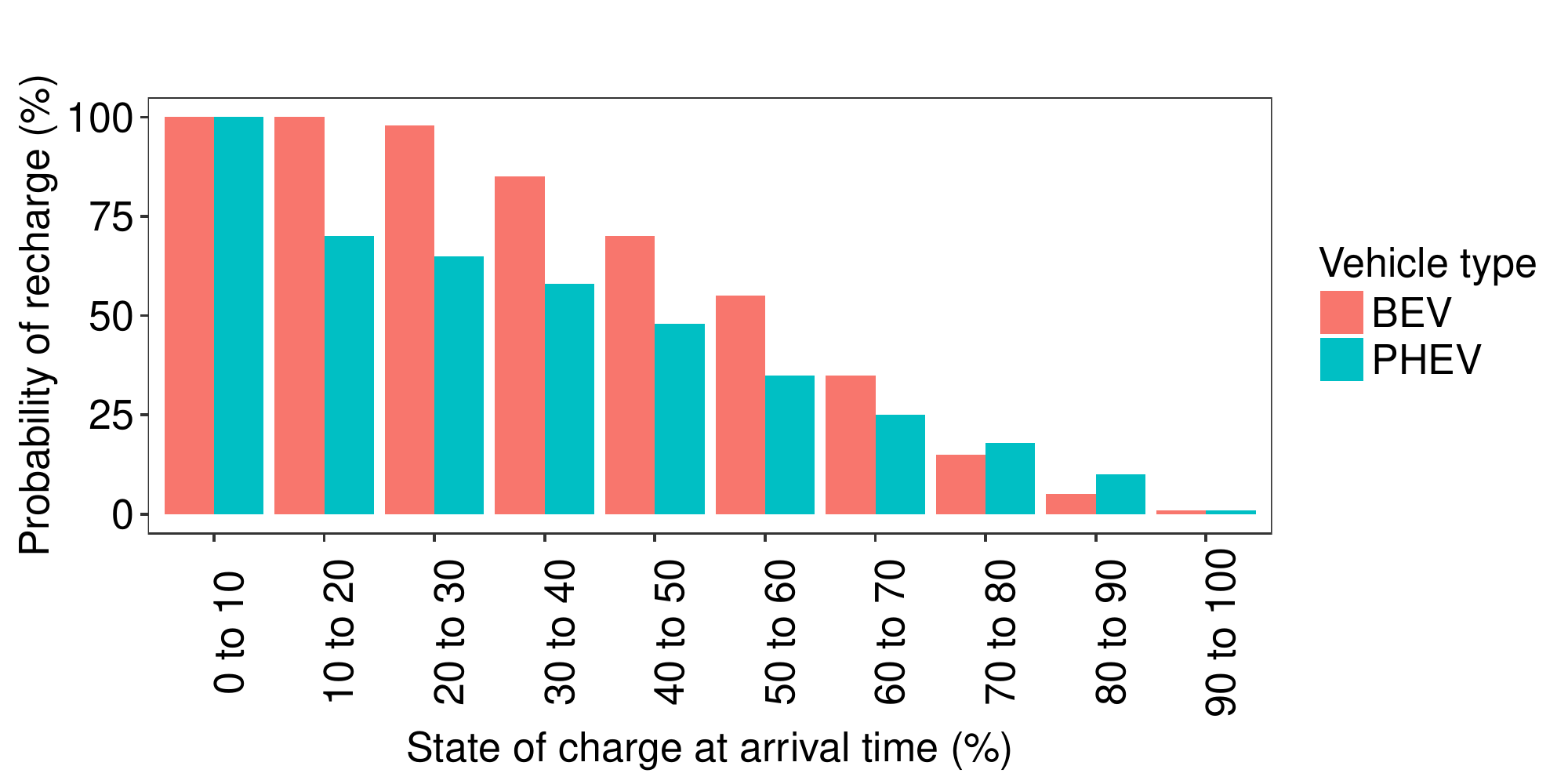}
	\caption{Probability of recharging as a function of the battery SOC at arrival time; Source: \cite{brooker2015identification}.}
	\label{SOC}
\end{figure}

A recent study analyzing two years of data from January 2011 to December 2013 of charging events that occurred away from home concluded that Nissan Leaf (pure battery electric vehicle, BEV) drivers prefer to charge their vehicles before their battery SOC drops to lower levels while Chevrolet Volt (a plug-in hybrid electric vehicle, PHEV) drivers tend to start recharging when there is a little charge in the battery since they rely on both electric motor and internal combustion engine \cite{brooker2015identification}. Fig. \ref{SOC} compares the probability of recharging for different values of battery SOC at the time of arrival for Nissan Leaf and Chevrolet Volt. 

\subsection{Dwell time}

\begin{figure}[!htbp]
	\centering
	\includegraphics[scale = 0.3]{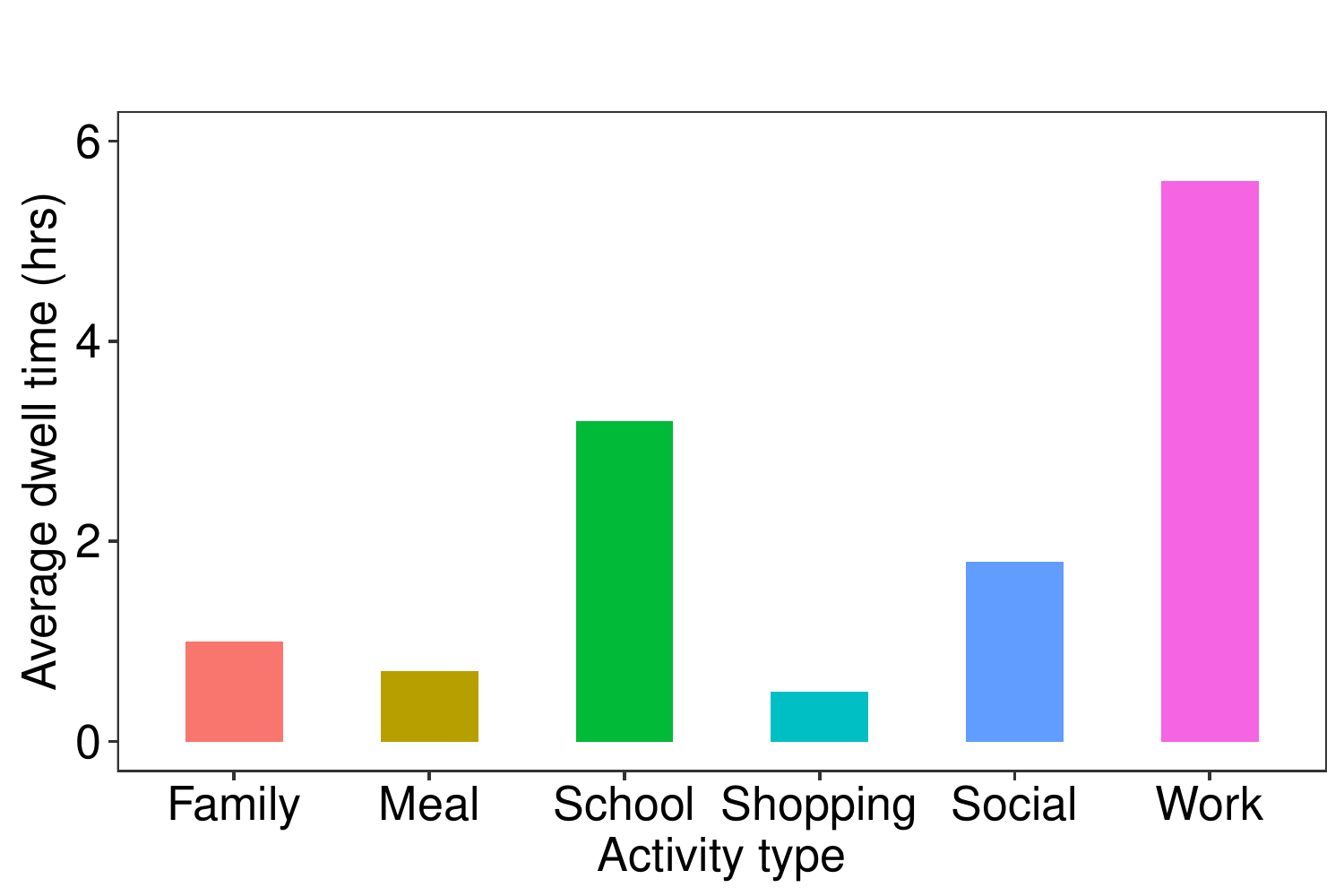}
	\caption {Average dwell time for activity types; Source: \cite{brooker2015identification}.}
	\label{Dwell}
\end{figure}

We define six different destination categories based on NHTS (National Household Travel Survey) data: Work, Social, Family, Meal, Study, and Shopping. Fig. \ref{Dwell} shows average time that people tend to park their vehicles based on their activity type \cite{krumm2012people}. Zhong et al. \cite{zhong2008studying} concluded that Weibull, log-normal and log-logistic distributions are the best distributions for modeling duration of weekday and weekend activities. While their analysis shows that model type and parameters or both might be different for an activity in weekday versus weekend, they found Weibull distribution the most applicable one. In addition, they found that certain activities such as social and shopping tend to last longer during weekends. Weibull distribution is used in our analysis to estimate parking duration of EV drivers considering average staying time, and we have also differentiated the durations of all weekday and weekend activities except meal activity.

\subsection{Weekday vs. weekend}

\begin{figure}[!htbp]
	\centering
	\includegraphics[scale = 0.3]{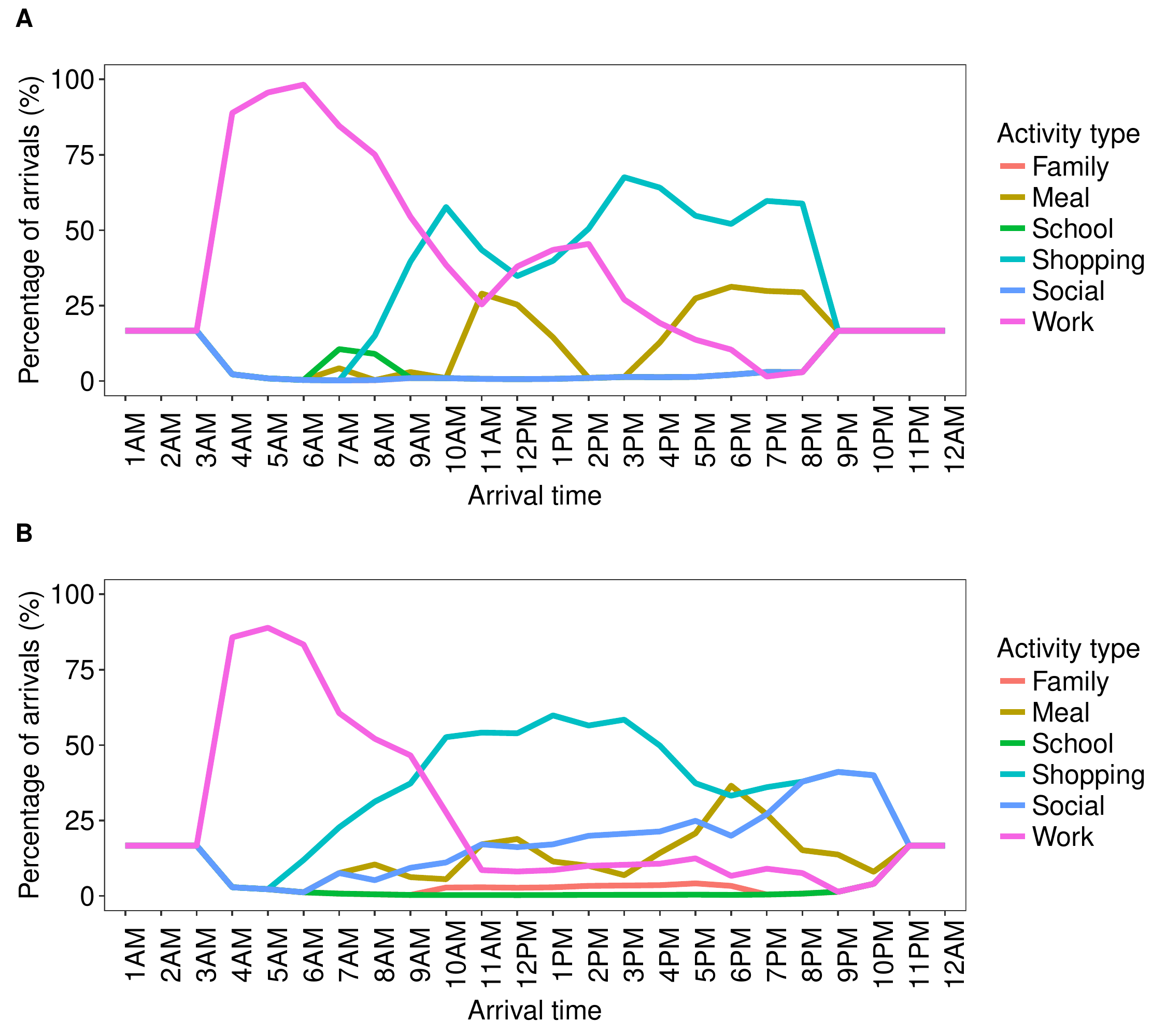}
	\caption {The expected breakdown of vehicle arrival percentages in A) weekdays and B) weekends; Sources: \cite{brooker2015identification} and \cite{krumm2012people}.}
	\label{Arrival}
\end{figure}

Demand pattern for public EV charging service can vary from day to day since people tend to attend social events, visit their families and go to shopping centers more during weekends than weekdays, in which demand mostly consists of people traveling to work or school. Fig. \ref{Arrival} confirms that demand for charging stations depends on time and type of day. During weekdays, maximum load occurs in morning when people are arriving at work or school while maximum demand usually happens around noon during weekends when people are going to shopping malls and social places. According to \cite{ozdemir2015distributed}, the best fitted distribution for arrival time to parking lot is a Weibull distribution. Hence, without loss of generality, we recommend the use of two Weibull distributions to estimate the arrival time of EVs to parking lots during weekdays and weekends.

\subsection{Preference for charging away from home}

Analysis by Idaho National Laboratory on data from 2012 and 2013 over 4,000 Leafs and 1,800 Volts across the U.S. shows that 13\% of Leaf drivers and 5\% of Volt drivers only charge their vehicles at home. This indicates that vast majority of drivers intend to use publicly accessible charging stations. This analysis also shows that although many people that drive more daily miles tend to charge their vehicles in places other than their homes, the effect of daily miles traveled on the chance of charging away from home is small. Hence, without loss of generality, we do not consider the effect of driving distance to community as a factor that affects the chance of using EV charging stations. 

\subsection{EV market penetration}

There are many social, environmental and economic factors that can significantly contribute to the increasing market share of different types of EVs. The survey in \cite{carley2013intent} about adult drivers in large U.S. cities in fall 2011 comprehended factors affecting the purchase of a plug-in EV. Besides demographic variables that can strongly predict intent of purchase, their results show that the presence of a charging station inside the community is the only awareness variable that has a significant effect on intent of purchase. Environmental Protection Agency estimated that 3.5\% of the vehicle fleet will be BEV or PHEV in the 2022-2025 time frame \cite{EPA2016midterm}.

\subsection{Willingness to walk}

The drivers' willingness to walk can be affected by their socio-demographic characteristics such as age, gender, education level and occupation. Many researchers have used distance decay function that shows the willingness to walk or bike as a distance towards different types of destinations. The parameter of this decay function depends on the activity type. Estimation results from \cite{yang2012walking} confirm that negative exponential distribution can better describe walking trips over short distances than other distributions such as Gaussian. They specify the distance decay function as 
\begin{equation}
	P(d) = e^{-\beta \times d}
\end{equation}
which shows the percentage of people willing to walk $d$ or longer distances than $d$. They used 2009 NHTS data to estimate the decay parameter $\beta$ for different groups and trip purposes. Their analysis shows that people are more willing to walk for recreation, social events and work activities rather than for studying, shopping or eating meal. Table \ref{WalkingPreference} shows the parameters of distance decay function influenced by variations in natural and built environment factors. The effects of season, region and community size on willingness to walk patterns are considered as well.

\begin{figure}[!htbp]
	\centering
	\includegraphics[scale = 0.3]{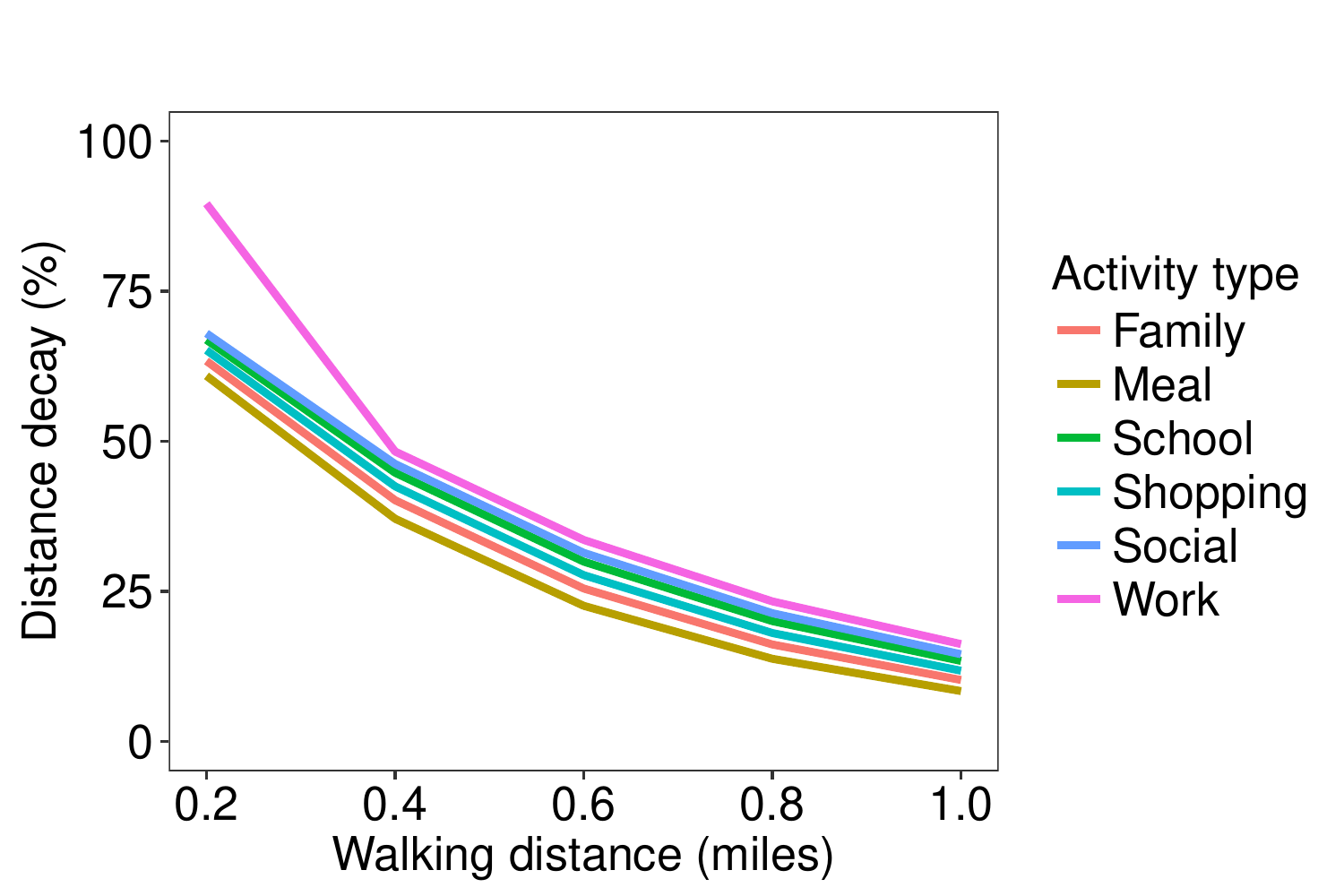}
	\caption {Distance decay function for walking trips to different destination types; Source: \cite{yang2012walking}.}
	\label{Walk}
\end{figure}

\begin{table}[!htbp]
	\caption{Estimated parameter for distance decay function}
	\centering
	\begin{tabular}{| c | c | c |} 
		\hline
		Factor & Category & $\beta$ \\ 
		\hline
		\multirow{4}{*}{Season} & Winter (Dec to Feb) & 1.88 \\
		& Spring (Mar to May) & 1.68 \\ 
		& Summer (Jun to Aug) & 1.64 \\
		& Autumn (Sep to Nov) & 1.7 \\ 
		\hline
		\multirow{4}{*}{Region} & Northeast & 1.85 \\
		& Midwest & 1.65 \\ 
		& South & 1.76 \\
		& West & 1.65 \\ 
		\hline
		\multirow{3}{*}{Community} & Town and country & 1.68 \\
		& Suburban & 1.63 \\
		& Urban and second city & 1.78\\  
		\hline
	\end{tabular}
	\label{WalkingPreference}
\end{table}

Our research aims at maximizing coverage of demand for public EV charging network in an urban/community area by proposing a two-stage stochastic programming model considering uncertainties in EV total flow, arrival and departure time, battery SOC at arrival time, preference for charging EV away from home, and walking preference patterns in the community.

\section{Model formulation and solution approach}\label{model}

\subsection{Model formulation}

Two-stage stochastic programming is a common approach for modeling problems that involve uncertainty in decision making. First-stage decision variables represent `here-and-now' decisions which are determined before the realization of randomness, and the second-stage decisions are determined after scenarios representing uncertainties are presented. In our model, binary variables in the first-stage determine the parking lots, and number of charging station installations for the selected parking lots. In the second-stage, a recourse decision is made on assigning EV drivers to one of their preferred parking lots based on their willingness to walk so that the expected access of EV drivers to public charging network is maximized.  

We first define the following model sets, parameters and variables:

\begin{itemize}
	\item Sets
	\begin{itemize}
		\item $S$: Set of parking lots, indexed by $s \in S$. 
		\item $L_{s}$: Set of number of charging stations in a parking lot $s$, indexed by $l \in L_{s}$.
		\item $B$: Set of buildings, indexed by $b \in B$.
		\item $T$: Set of time slots, indexed by $t \in T$.
		\item $\Gamma$: Set of arrival and departure times, indexed by $\gamma (t) \in \Gamma$ containing time slot $t\in T$.
		\item $\Omega$: Set of scenarios, indexed by $\omega \in \Omega$.
	\end{itemize} 
	\item Model Parameters
	\begin{itemize}
		\item $p$: Number of parking lots to be considered for installing charging stations.
		\item $m_l$: Number of charging stations, $l \in L_s$.
		\item $d_{\gamma (t),b,s}(\omega)$: Demand with arrival and departure time set of $\gamma (t) \in \Gamma$ for a given $t \in T$ for a building $b$ that are willing to park their vehicle in parking lot $s \in S', S' \subset S$ in a scenario $\omega \in \Omega$.
	\end{itemize}
	\item First Stage Decision Variables
	\begin{itemize}
		\item $x_s$: 1, if parking lot $s \in S$ is considered for installing charging stations; 0, otherwise.
		\item $z_{l,s}$: 1, if $l \in L_s$ charging stations are installed in parking lot $s \in S$.
	\end{itemize}
	\item Second Stage Decision Variables
	\begin{itemize}
		\item $y_{\gamma (t),b,s}(\omega)$: Proportion of demand with arrival and departure time set of $\gamma (t) \in \Gamma$ for a building $b$  willing to charge their vehicle in parking lot $s \in S', S' \subset S$ in a scenario $\omega \in \Omega$.
	\end{itemize}
\end{itemize}

The two-stage stochastic programming model is as follows:
\begin{alignat}{3}
& \text{First-Stage Model:} \label{First-Stage Model} && \nonumber \\
& \text{ Max } E_{\Omega}[\varphi(x,z,\omega)]  && \\ \label{fs-obj}
& \text{s.t. } && \nonumber \\
& \sum_{s\in S} x_{s} = p && \\ \label{fs-obj-eq3}
& z_{l,s} \leq x_{s} \hspace{4.65cm} \forall  s\in S, l\in L_{s} && \\  \label{fs-obj-eq4}
& \sum_{l\in L_{s}} z_{l,s} \leq 1 \hspace{5cm} \forall s \in S &&  \\  
& x_{s}, z_{l,s} \in \{0,1\} \hspace{3.6cm} \forall s\in S, l\in L_{s} \label{fs-obj-eq6} &&  \\
& \text{Second-Stage Model:} && \nonumber \\
& \text{ Max } \varphi(x,z,\omega) && \nonumber \\
& = \sum_{t \in T} \sum_{\gamma (t) \in \Gamma} \sum_{b \in B} \sum_{s \in S} y_{\gamma (t),b,s}(\omega)d_{\gamma (t),b,s}(\omega) \label{Second-Stage Model} && \\
& \sum_{\gamma (t) \in \Gamma} \sum_{b \in B} y_{\gamma (t),b,s}(\omega)d_{\gamma (t),b,s}(\omega) \leq \sum_{l \in L_{s}} m_{l}z_{l,s} && \nonumber \\
& \hspace{5.8cm} \forall s \in S, t \in T && \\
& \sum_{s \in S} y_{\gamma (t),b,s}(\omega) \leq 1 \hspace{1.5cm} \forall t \in T, \gamma (t) \in \Gamma, b \in B && \\
& 0 \leq y_{\gamma (t),b,s}(\omega) \leq 1 \hspace{0.45cm} \forall t \in T, \gamma (t) \in \Gamma, b \in B, s \in S &&
\end{alignat}

In this model, first-stage decisions are made regarding the locations of charging stations and charging capacity in each location. The first-stage objective function (\ref{First-Stage Model}) maximizes the expected access, and $E_{\Omega}$ is an expectation operator, and $E_{\Omega}[\varphi(x,z,\omega)]$ represents $\sum_{\omega \in \Omega}p_\omega \varphi(x,z,\omega)$, where $p_\omega$ is probability of occurrence for scenario $\omega$, and $\sum_{\omega \in \Omega}p_{\omega}=1$. Constraint (3) ensures that $p$ parking lots are selected to install EV charging stations. Constraints (4) and (5) determine charging capacity in any parking lot that is selected for providing EV charging service. Constraints (6) define the feasible set for the binary first-stage variables. In the second-stage, recourse decisions are made to maximize the coverage of potential EV traffic flows based on the decisions chosen in the first-stage and a realization $\omega \in \Omega$. Constraints (8) describe the supply-demand balance restrictions. They ensure that demand that has arrival and departure time set of $\gamma (t)$ and are assigned to parking lot $s$ for EV charging does not exceed the charging capacity in parking lot $s$. Constraints (9) state that demand with arrival and departure time set of $\gamma (t)$ can be assigned to at most one parking lot for EV charging. Constraint set (10) are the non-negativity constraints. Though we have not considered any budgetary restrictions, such constraints can be added to the first-stage model if appropriate.

\subsection{Solution approach}

\subsubsection{Sample average approximation}

According to \cite{mak1999monte}, unless there are small number of scenarios that can represent uncertainties in a problem, it is usually impossible to solve a stochastic programming problem. They showed that optimal solution of stochastic programming can be approximated by a sample of scenarios much smaller than the actual size of scenarios and this approximation monotonically improves as we increase the number of scenarios. SAA is also an effective approach when sufficient number of scenarios to estimate optimal solution is unknown. SAA was proposed by \cite{mak1999monte} and for the sake of completeness, we provide the procedure for sample average approximation method as follows:

\begin{enumerate}
	\item Estimating an upper bound for the optimal solution:
	\begin{itemize}
		\item Generate $M$ independent sample sets of scenarios each of size $N$, i.e., ($\omega^1_j,\omega^2_j,...,\omega^N_j$) for $j = 1,2,...,M$
		\item For each sample set $j = 1,2,...,M$, find the optimal solution:
		\begin{equation}
		v^j_{N} = \frac{1}{N} \sum_{i=1}^{N} \varphi(x,z,\omega^i_j).
		\end{equation}
		\item Compute the followings:
		\begin{equation}
		\overline v_{N,M} = \frac{1}{M} \sum_{j=1}^{M} v^j_{N}
		\end{equation}
		\begin{equation}
		\sigma^2_{\overline v_{N,M}} = \frac{1}{M(M-1)} \sum_{j=1}^{M} (v^j_{N} - \overline v_{N,M})^2.
		\end{equation}
		\\
		The expected value of $v_N$ is greater than or equal to the optimal value $v^*$. Since the sample average $\overline v_{N,M}$ is an unbiased estimation of the expected value of $v_N$, $\overline v_{N,M}$ provides an upper statistical bound for the optimal solution. 
	\end{itemize}
	\item Estimating a lower bound for the optimal solution:
	\begin{itemize}
		\item If ($\overline x,\overline z$) is a feasible solution for the first-stage problem, then $f(\overline x,\overline z) \leq v^*$. Hence, choosing any feasible solution of the first-stage problem will provide a lower statistical bound for the optimal value. 
		\item Choose a sample of size $N'$ of scenarios, much larger than $N$, i.e., ($\omega^1,\omega^2,...,\omega^{N'}$) and independent of samples to find the upper limit and estimate the objective function:
		\begin{equation}
		f(\overline x, \overline z) = \frac{1}{N'} \sum_{i=1}^{N'} \varphi(x,z,\omega^i) 
		\end{equation}
		\item Compute the variance for this estimation:
		\begin{eqnarray}
		\sigma^2_{N'}(\overline x,\overline z) = \frac{1}{N'(N'-1)} \sum_{i=1}^{N'} (\varphi(x,z,\omega^i) \nonumber \\
		- f(\overline x, \overline z))^2.
		\end{eqnarray}
	\end{itemize}
	\item Estimating the optimality gap: 
	\begin{itemize}
		\item  Use the upper bound and the lower bound that are computed in previous steps to estimate the optimality gap:
		\begin{equation}
		gap_{M,N,N'}(\overline x,\overline z) = \overline v_{N,M} - f(\overline x,\overline z).
		\end{equation}
	\end{itemize}
	\item Checking the quality of the estimated optimality gap:
	\begin{itemize}
		\item Variance of the estimated optimality gap can be found by
		\begin{equation}
		\sigma^2_{gap} = \sigma^2_{\overline v_{N,M}} + \sigma^2_{N'}(\overline x,\overline z)
		\end{equation}
	\end{itemize}
\end{enumerate}

\subsubsection{Heuristic}

SAA requires high computational resources, hence we developed a heuristic to solve large-scale problems efficiently. This heuristic is inspired by a score measure introduced by \cite{tsiligirides1984heuristic}. The score incorporates charging capacity of each parking lot as well as its distance to other parking lots. The heuristic consists of a construction phase during which we build an initial solution, and an improvement phase where we employ local search moves to find a better solution. The pseudo-code of the heuristic is presented as follows: 

\begin{algorithm}[!htbp]
	\caption{Pseudo-code of the heuristic}
	\begin{algorithmic}[1]
		\State $bestsolution \gets \emptyset$.
		\For {$s\gets 1$ to $NumberofParkingLots$}:
		\State Compute score measure $r_{s}$.
		\EndFor
		\State \textbf{Construction phase}:
		\State $initialsolution \gets \emptyset$
		\State Compute attractiveness ratio $\rho_{s}$ for all parking lots.	
		\State Add parking lots to the initial solution in decreasing order of the attractiveness ratio until $p$ parking lots are selected.
		\State \textbf{Improvement phase}:
		\State $currentsolution \gets initialsolution$	 
		\While{$f(currentsolution)$ can be improved}
		\State remove-insert($currentsolution$)
		\EndWhile\label{euclidendwhile}
		\State Store best solution found so far.
	\end{algorithmic}
\end{algorithm}

In the construction phase, a score measure for each parking lot as a potential location for installing charging stations is calculated as:
\begin{equation}
	r_s = \sum_{s,s' \in S, s' \neq s} c_s e^{-\beta d_{s,s'}}
\end{equation}
where $\beta$ is a user parameter. The score is measured as an incentive for the charging capacity ($c_s$) of each parking lot, and distance ($d_{s,s'}$) to other parking lots as a cost. If a parking lot has more capacity for installing charging stations and is nearer to other parking lots, its score would be higher. 

To consider randomness in constructing the initial solution, we use a set of sample scenarios to get the probability of parking lot $s$ being chosen as one of the optimal locations for installing charging stations. This estimated probability for parking lot $s$, $q_s$, is computed based on the fraction of scenarios in which parking lot $s$ is among the optimal locations. The attractiveness measure of parking lot $s$, $\rho_s$, is computed by multiplying this probability to the corresponding score measure: 
\begin{equation}
	\rho_s = r_s q_s
\end{equation}

Parking lots will be added to the initial solution in a decreasing order of attractiveness measure until $p$ parking lots are selected. In the improvement step, we use local search method of remove-insert procedure. For every parking lot that is already in the initial solution, we replace it with one of the parking lots that has not been selected based on a parking lot that has the highest attractiveness measure. This process is continued until there is no improvement in the objective function. We repeat this procedure for all parking lots that are selected in the initial solution and store the best value found for the objective function.   

\section{Case study and computational experiments}\label{case}

To demonstrate the efficacy of the proposed approach, our case study investigates the community area data of Detroit midtown area in Michigan, U.S. There is a wide range of employment types (type of final destinations) in this area, it attracts a lot of traffic, and is characterized by an urban university, commercial offices, hospitals, and museums. This area includes 135 buildings among which 67 are office buildings, 12 are social places, 5 are family related buildings, 4 are restaurants, 44 are schools buildings and 3 are shopping places. There are 32 parking lots that are considered as potential locations for installing EV charging stations. We assume that parking lots are open between 6am and 6pm, and have different capacities for installing charging stations. The center of each parking lot is considered as our candidate for installing a charging station, and Euclidean distance is used to measure distance between any two points in the community. Data from Southeast Michigan Council of Governments shows that average annual daily traffic of Detroit midtown area is approximately between 10,000 and 20,000 and like \cite{capar2013arc}, we assume that total daily traffic of this community follows a uniform probability distribution.

\begin{figure}[!htbp]
	\centering
	\includegraphics[scale = 0.25]{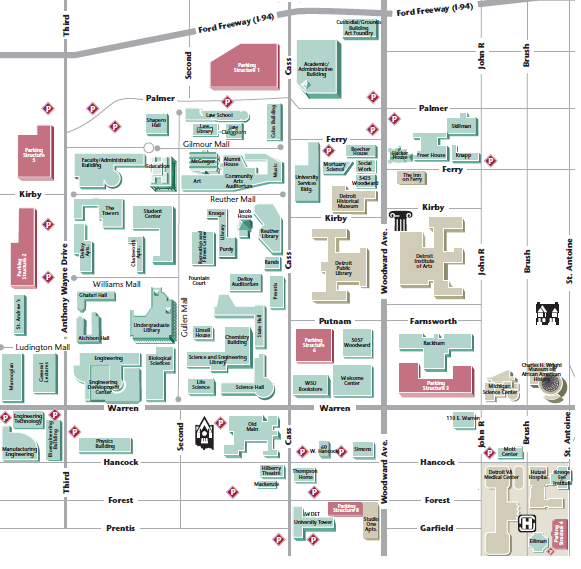}
	\tiny \caption {Part of Detroit midtown area used for our analysis.}
\end{figure}

Based on Environmental Protection Agency (EPA) analysis, we examine cases in which EVs constitute 3\% and 5\% of the light-duty vehicle fleet. According to \cite{vergis2015comparison}, weather/climate is positively correlated with BEV market share. Since our case study is done in an area with low winter temperatures, BEV market share is considered lower than PHEV market share. Two cases are constructed for our computational experiments. In the first case, we assume that the market share is 1\% and 2\% for BEVs and PHEVs, respectively. In the second case, these market shares are assumed to be 2\% and 3\%, respectively.

In this study, negative exponential distribution functions estimated by \cite{yang2012walking} are used to describe willingness to walk patterns for various activity types which considers the effects of season and community size in U.S. Drivers are randomly assigned to a parking lot that is within their walking distance preference. In both cases of EV market share, 13\% of total demand is not considered in our model since there is no parking lot within their walking distance preference, and also it is difficult to track their walking distance to their final destination if they use other EV charging sources placed in streets, etc.  Four different values (2,4,6 and 8) for number of parking lots ($p$) to install charging stations are considered. The optimization models for SAA and heuristic were implemented in Python 2.7 using Gurobi 6.5.1 software for solving optimization problems. All the computations were performed using a system with Intel (R) Xeon(R) CPU 3.10 GHz and 24GB RAM.

\subsection{Scenario construction}

For the two-stage model, uncertainties are modeled by use case scenarios. A scenario represents a single day of public EV charging service and is influenced by short-term (weekday vs. weekend) and long-term (seasonal) variations, and total number of EVs arriving to the community. The probability of occurrence for a scenario is based on a uniform probability distribution. Without loss of generality, we assume that any given scenario day can belong to winter, spring, summer, and autumn seasons with equal probability. 

In each scenario, a random number from $U(0,1)$ determines type of each vehicle in a community, and if the random number is less than BEV market share, between BEV market share and sum of BEV and PHEV market shares, or greater than sum of BEV and PHEV market shares, then the vehicle is assumed to be a BEV, a PHEV or an ICE (internal combustion engine), respectively. If it is an EV, Weibull distributions with parameters (8,3) and (13,4) are used to determine arrival time of EV drivers to the community in a weekday and weekend day, respectively. As explained in the earlier section, the purpose of arrival for a driver is determined based on arrival time and distributions. Furthermore, Weibull distribution is used to estimate duration of various types of weekday and weekend activities. Table \ref{Params} represents the parameters for this distribution based on type of activity. In this table, the first and second numbers represent the shape parameter and the scale parameter, respectively.

 \begin{table*}
 	\caption{Weibull distribution parameters for activity duration}
 	\centering
 	\begin{tabular}{| c | c | c | c | c | c | c |}  
 		\hline
 		Type of day & Work & Social & Family & Meal & School & Shopping \\ 
 		\hline
		Weekday & (5.89,10) & (1.89,10) & (1.05,10) & (0.79,2) & (3.61,2) & (0.56,2) \\
 		\hline
 		Weekend & (6.04,6) & (2.03,2) & (1.13,2) & (0.79,2) & (3.36,10) & (0.25,0.5) \\
 		\hline
 	\end{tabular}
 	\label{Params}
 \end{table*} 

For a final destination, each EV driver is randomly assigned to a target destination/building using a uniform distribution based on a driver's purpose of arrival to the community. A random number is generated from exponential distribution as shown in Fig. \ref{Walk} to determine each EV driver's willingness to walk distance based on his/her purpose of arrival, and also, community size and type of region are considered in willingness to walk distributions. If there is no parking lot within a driver's willingness to walk distance, then this demand is not considered in our model. In order to incorporate charging preference of EV drivers, uniform distribution $U(0,1)$ is used. If the random number is greater than 13\% for BEV or 5\% for PHEV, a driver's willingness to charge away from home is decided. Consistent with recommendations from \cite{dai2011optimization} and \cite{luo2013optimal}, without loss of generality, we assume that the initial battery SOC for vehicles arriving at the charging stations follows a normal distribution $N(0.3,0.1)$ with a mean 0.3 and standard variation 0.1. Based on battery SOC at arrival time, uniform distribution $U(0,1)$ is used to determine each EV driver's willing to charge EV at public charging stations. This is further compared with associated probability of recharge based on type of EV discussed in SOC section earlier. If the random number is less than or equal to the probability of recharge, that EV is considered as demand for EV charging network in the community. Similarly, multiple scenarios are constructed for the two-stage stochastic programming model to simulate the arrival pattern, battery SOC, dwell time, charging preference and willingness to walk in the community.

\subsection{SAA settings}

To estimate an upper bound for expected accessibility to public EV charging stations, $N$ = 30, 50 and 100 scenarios are used and this is repeated $M$ = 20 times. The average of these 20 runs is an estimate of upper bound on the accessibility. A sample of $N'$ = 1,000 scenarios, which are separate from those that were used to get the upper bound, is used to estimate a lower bound for the optimal solution. Computation times for each test problem along with the heuristic performance are summarized in Tables \ref{SAA12} and \ref{SAA23}. The computation times show that the optimization model using SAA method is able to solve problems with eight optimal locations in less than five hours. In these tables, UB (\%) and LB (\%) represent upper and lower bounds for expected accessibility to public EV charging service using SAA method. Gap (\%) and gap SD indicate the differences between upper and lower bounds and standard deviation, respectively. Opt(s) is the running time of SAA. The best solution found by our heuristic for upper bound of the objective function and its running time are shown as Heuristic (\%) and Heuristic (s).

\begin{table*}[!htbp]
	\caption{SAA performance when $(M,N')$ = (20,1,000) and (BEV,PHEV) = (1\%,2\%)}
	\centering
	\begin{tabular}{|c c c c c c c c c|} 
		\hline
		$p$ & $N$ & UB (\%) & LB (\%) & gap (\%) & gap SD & Opt (s) & Heuristic (\%) & Heuristic (s) \\
		\hline
		\multirow{3}{1em}{2} & 30 & 57.98 & 56.59 & 2.39 & 0.0064 & 397 & 57.98 & 68 \\
		& 50 & 58.70 & 58.25 & 0.77 & 0.0062 & 1,226 & 58.70 & 74 \\
		& 100 & 58.56 & 58.54 & 0.02 & 0.0055 & 4,564 & 58.56 & 93 \\
		\hline
		\multirow{3}{1em}{4} & 30 & 73.89 & 73.42 &	0.63 & 0.0056 & 720 & 73.88 & 114 \\
		& 50 & 74.61 & 73.85 & 1.02 & 0.0041 & 1,759 & 74.61 & 131 \\
		& 100 & 74.59 & 73.74 &	1.14 & 0.0040 & 7,406 & 74.59 & 193 \\
		\hline
		\multirow{3}{1em}{6} & 30 & 83.97 & 83.62 &	0.35 & 0.0039 & 1,071 & 83.21 & 160 \\
		& 50 & 84.11 & 83.80 & 0.31 & 0.0034 & 2,173 & 83.17 & 186 \\
		& 100 & 83.40 & 83.30 & 0.10 & 0.0031 & 9,572 & 82.86 & 303 \\
		\hline
		\multirow{3}{1em}{8} & 30 & 91.16 & 90.61 & 0.61 & 0.0026 & 1,124 & 90.28 & 185 \\
		& 50 & 91.13 & 90.78 & 0.38 & 0.0021 & 3,099 & 90.18 & 245 \\
		& 100 & 90.87 & 90.86 &	0.02 & 0.0018 & 12,832 & 90.11 & 414 \\
		\hline
	\end{tabular}
	\label{SAA12}
\end{table*}

\begin{table*}[!htbp] 
	\caption{SAA performance when $(M,N')$ = (20,1,000) and (BEV,PHEV) = (2\%,3\%)}
	\centering
	\begin{tabular}{|c c c c c c c c c|} 
		\hline
		$p$ & $N$ & UB (\%) & LB (\%) & gap (\%) & gap SD & Opt (s) & Heuristic (\%) & Heuristic (s) \\
		\hline
		\multirow{3}{1em}{2} & 30 & 50.42 & 50.00 & 0.85 & 0.0056 & 462 & 50.42 & 82 \\
		& 50 & 50.91 & 50.10 & 1.58 & 0.0054 & 1,141 & 50.91 & 87 \\
		& 100 & 50.91 & 50.31 &	1.17 & 0.0048 & 4,761 & 50.91 & 106 \\
		\hline
		\multirow{3}{1em}{4} & 30 & 63.35 & 63.16 &	0.30 & 0.0064 & 1,595 & 63.33 & 169 \\
		& 50 & 63.19 & 63.11 & 0.13 & 0.0063 & 3,644 & 63.19 & 211 \\
		& 100 & 63.46 & 63.42 & 0.07 & 0.0057 & 16,656 & 63.41 & 317 \\
		\hline
		\multirow{3}{1em}{6} & 30 & 72.56 & 71.55 & 1.39 & 0.0071 & 1,663 & 72.34 & 208 \\
		& 50 & 72.04 & 71.46 & 0.81 & 0.0059 & 3,246 & 71.84 & 273 \\
		& 100 & 71.82 & 71.40 & 0.58 & 0.0050 & 12,165 & 71.73 & 474 \\
		\hline
		\multirow{3}{1em}{8} & 30 & 78.91 & 78.49 &	0.52 & 0.0048 & 1,494 & 78.53 & 273 \\
		& 50 & 79.44 & 78.92 & 0.66 & 0.0045 & 2,908 & 79.01 & 374 \\
		& 100 & 79.12 & 78.69 &	0.54 & 0.0044 & 12,248 & 78.70 & 667 \\
		\hline
	\end{tabular}
	\label{SAA23}
\end{table*}

\subsection{Performance measures} 

Number of public chargers per capita could have a significant effect on both BEV market share and PHEV market share. In terms of monetary benefits for EV consumers, the average of total benefits across 25 major metropolitan areas is around \$2,800 per BEV and \$1,600 per PHEV \cite{lutsey2015assessment}. In order to deal with uncertainties in demand for public EV charging service and simulate the expected output measures with different number of chargers in the community, a set of 50 scenarios are generated and used for our analysis. We study two different cases for willingness to walk pattern in the community to generate optimistic and pessimistic bounds for level of walking in people that have access to public EV charging network. In the optimistic case, we assume that people are willing to walk long distances and will always choose the farthest available charging station to their final destination whereas in the pessimistic case people are willing to walk short distances and always choose the nearest available station to their building. Five different indicators are used to measure the performance of public EV charging placement: accessibility, lost demand, charging utilization, total walking distance, and walking distance per capita. Access is defined as the percentage of EV drivers that could charge their vehicles in public charging stations in the community, and lost demand is the percentage of EV drivers that are willing to use public EV network but there is not enough capacity to serve them. Charging utilization is the percentage of time that a charging station is being used by an EV. To assess walking patterns among people before network and after installing pubic EV charging stations, we use total walking distance and walking distance per capita measures.   

\begin{figure}[!htbp]
	\centering
	\includegraphics[scale = 0.3]{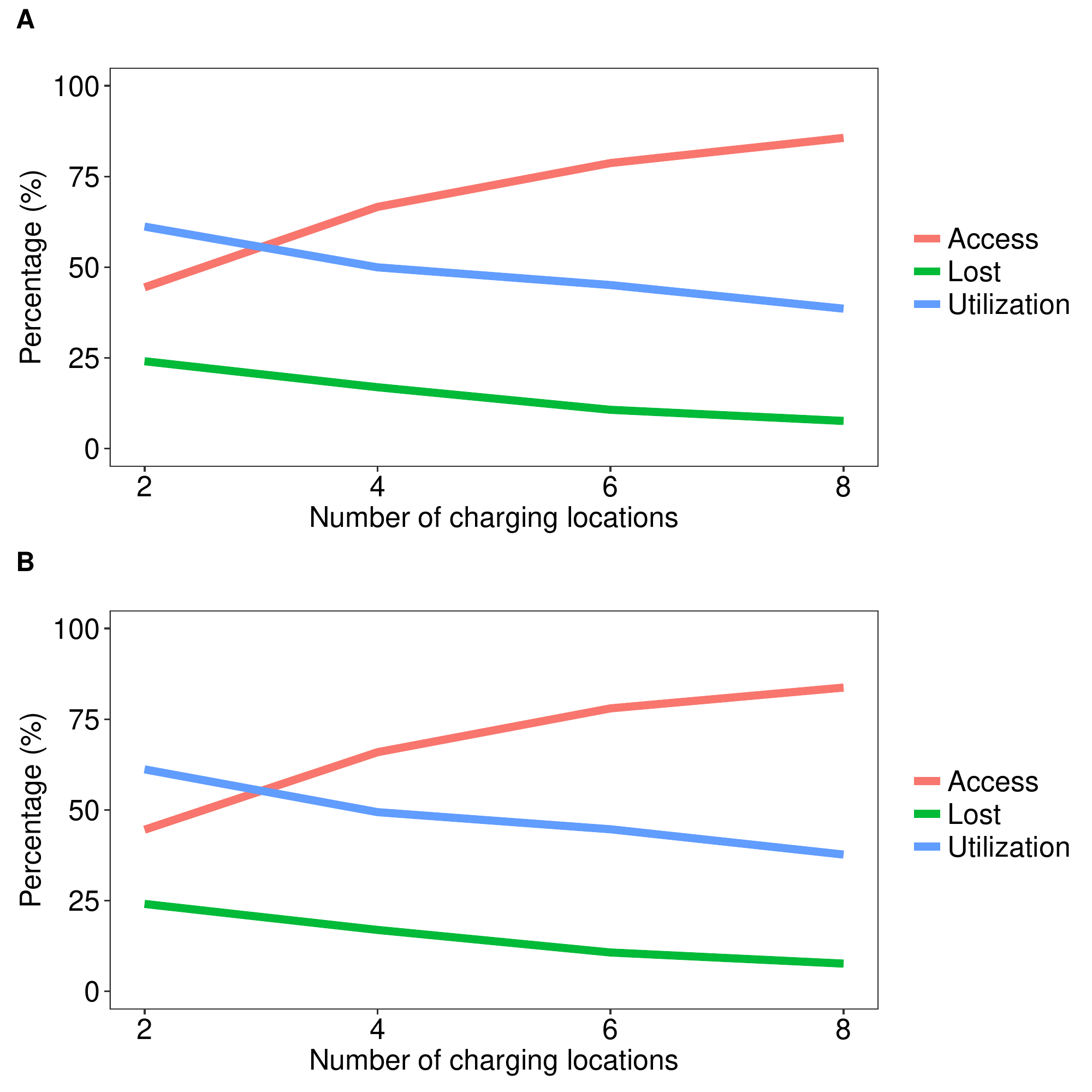}
	\caption {Percentage of accessibility, lost demand and charging utilization in A) optimistic and B) pessimistic cases when (BEV,PHEV) market shares are (1\%,2\%).}
	\label{Output12}
\end{figure}

\begin{figure}[!htbp]
	\centering
	\includegraphics[scale = 0.3]{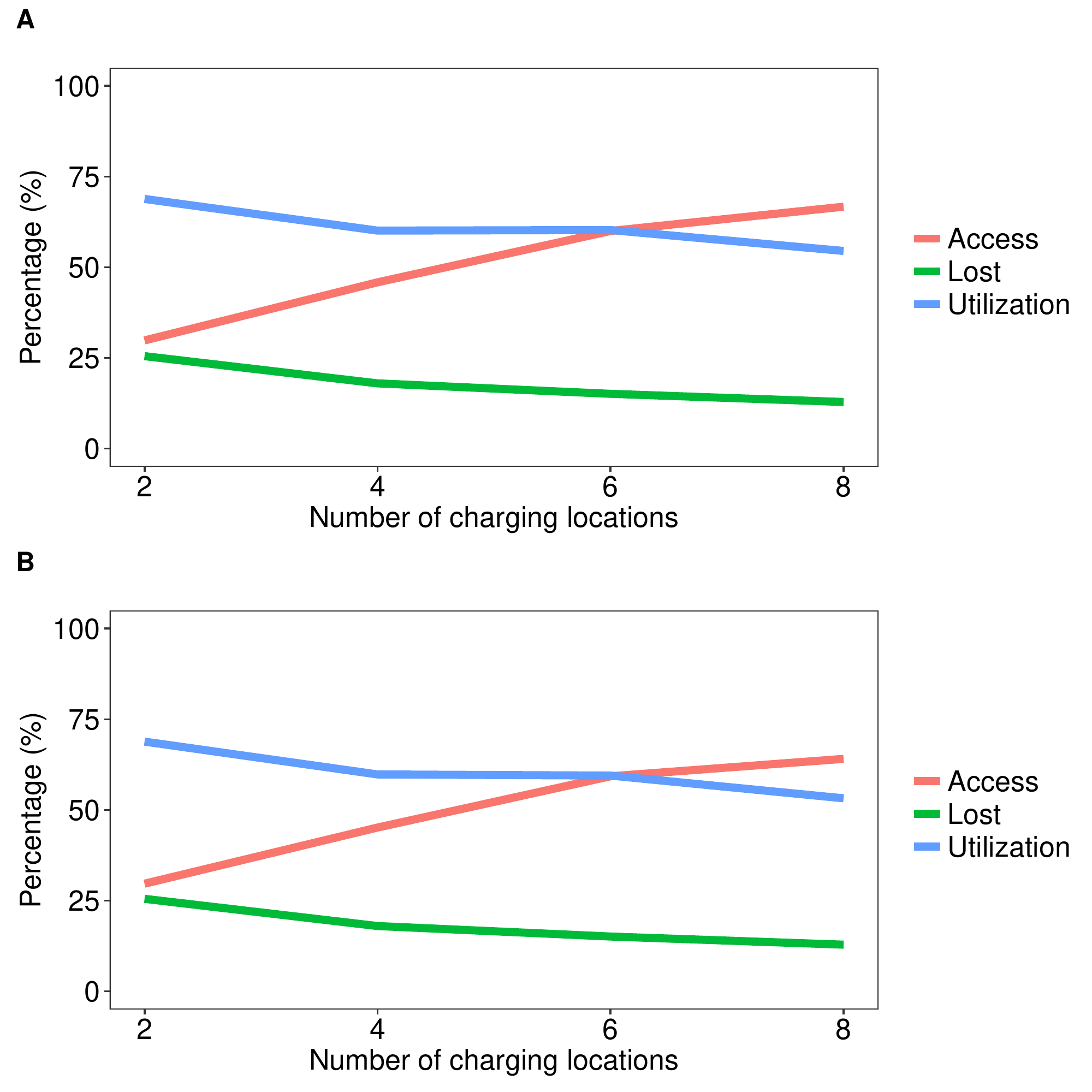}
	\caption {Percentage of accessibility, lost demand and charging utilization in A) optimistic and B) pessimistic cases when (BEV,PHEV) market shares are (2\%,3\%).}
	\label{Output23}
\end{figure}

As shown in Figs. \ref{Output12} and \ref{Output23}, accessibility to public charging service increases in both cases of EV market share as more charging stations are installed in the community but utilization level of these stations reduces simultaneously. Increase in EV market share can reduce accessibility to public charging network up to 32\% in both optimistic and pessimistic cases. However, this increase in demand will increase utilization level up to 41\% and lost demand up to 68\%.   

\begin{figure}[!htbp]
	\centering
	\includegraphics[scale = 0.3]{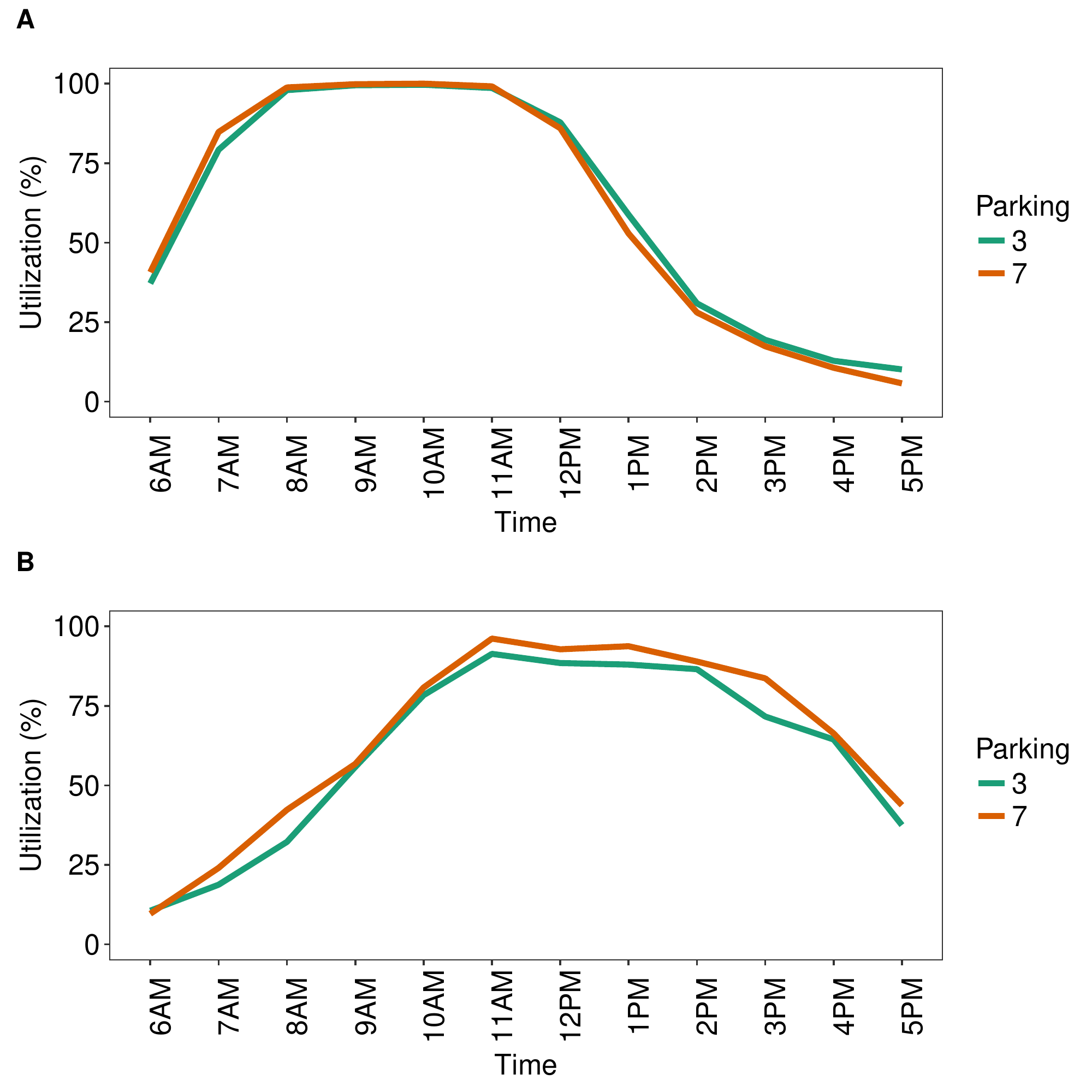}
	\caption {Average percentage of hourly utilization in A) weekdays and B) weekends in an optimistic case when $p$ = 2 and (BEV,PHEV) market shares are (1\%,2\%).}
	\label{HourlyUtilization212}
\end{figure}

\begin{figure}[!htbp]
	\centering
	\includegraphics[scale = 0.3]{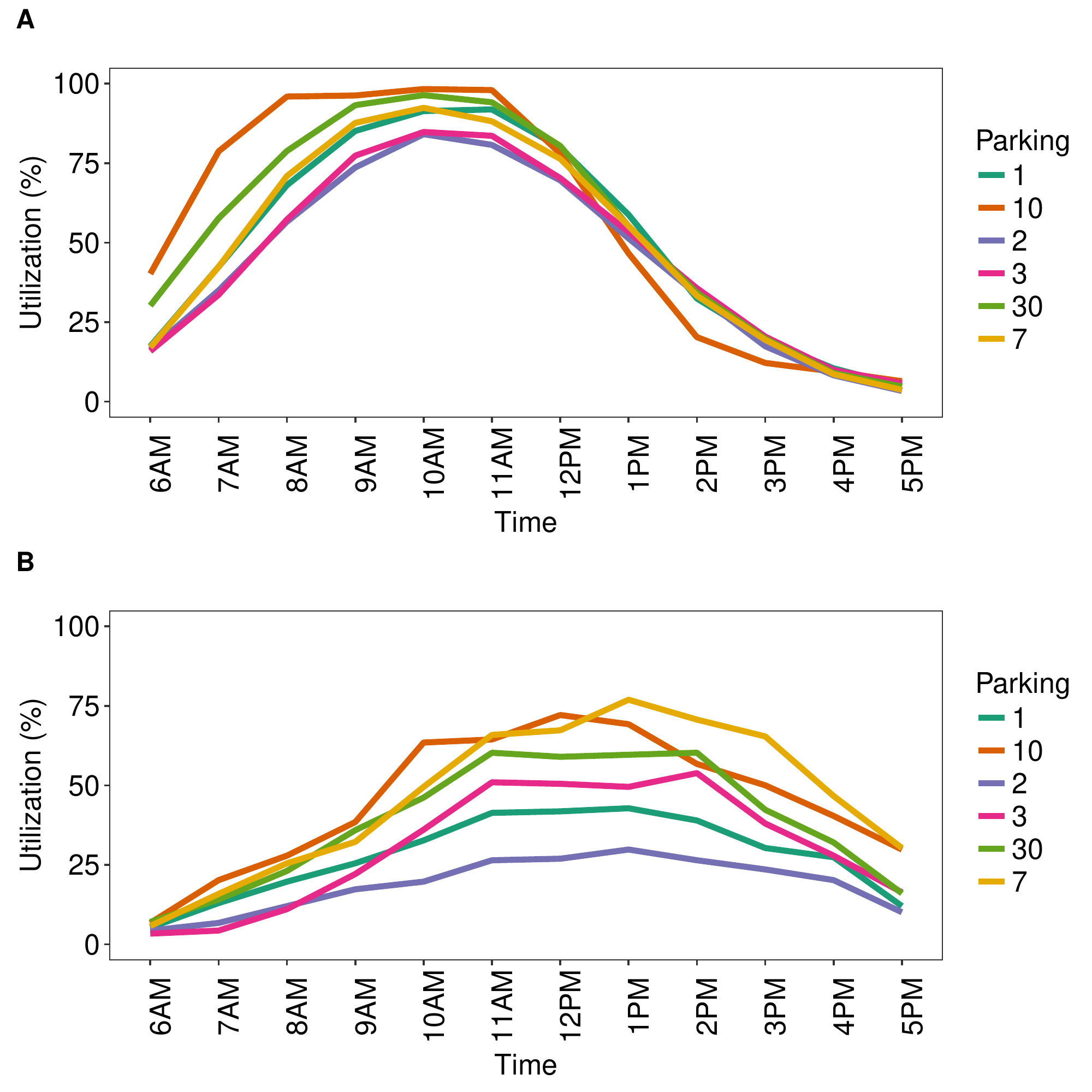}
	\caption {Average percentage of hourly utilization in A) weekdays and B) weekends in an optimistic case when $p$ = 6 and (BEV,PHEV) market shares are (1\%,2\%).}
	\label{HourlyUtilization612}
\end{figure}

Figs. \ref{HourlyUtilization212} and \ref{HourlyUtilization612} compare the average percentage of hourly utilization level of charging stations in weekdays versus weekends in an optimistic case of willingness to walk, and indicate a difference in utilization pattern from weekday to weekend. Utilization peaks around 8 in the morning during a weekday while it is around noon during weekend. These patterns match the expected arrival pattern of people to the community in weekdays and weekends. These plots also indicate that charging stations would not be fully utilized as more charging stations are available for EV drivers. This is important from revenue perspective since utilization level is among the major drivers of profitability of investment on public EV charging stations \cite{chang2012financial}.

An important measure of livability analysis via transportation is increasing the travel options so that people can meet at least a part of their travel needs through walking and biking, and improve their health condition \cite{marshall2013evaluation}. It has been estimated that a shift from driving to walking can save the average approximately 25c per vehicle-mile traveled and 50c under urban-peak condition, when emission and parking costs are high, in external costs such as traffic congestion, noise and air pollution \cite{litman2004economic}. Design of an effective EV charging network can also provide opportunities for people in a community to increase their level of physical activity.

\begin{figure}[!htbp]
	\centering
	\includegraphics[scale = 0.3]{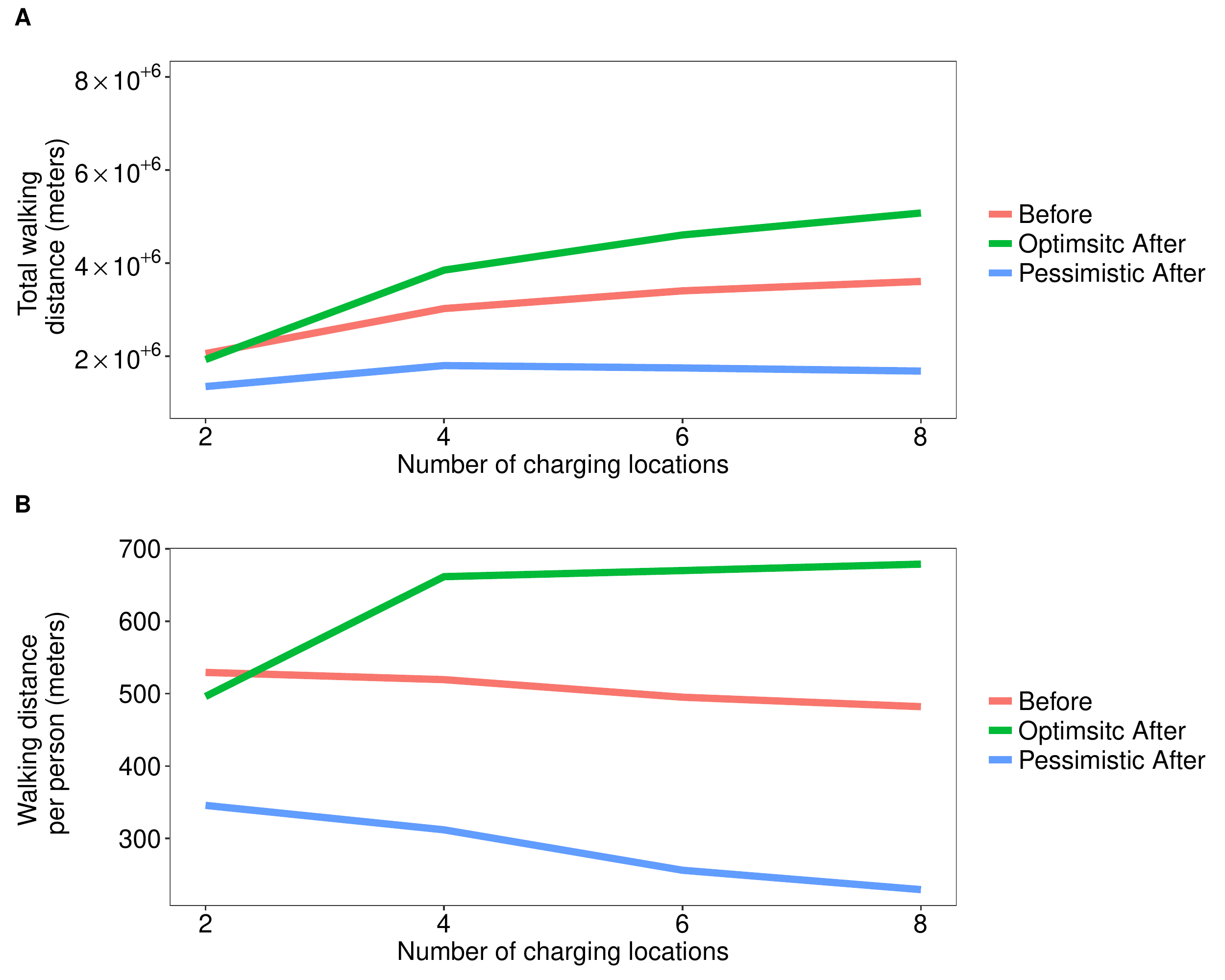}
	\caption {A) Total walking distance and B) walking distance per capita for people that have access to public EV charging service when (BEV,PHEV) market shares are (1\%,2\%).}
	\label{Walk12}
\end{figure}

\begin{figure}[!htbp]
	\centering
	\includegraphics[scale = 0.3]{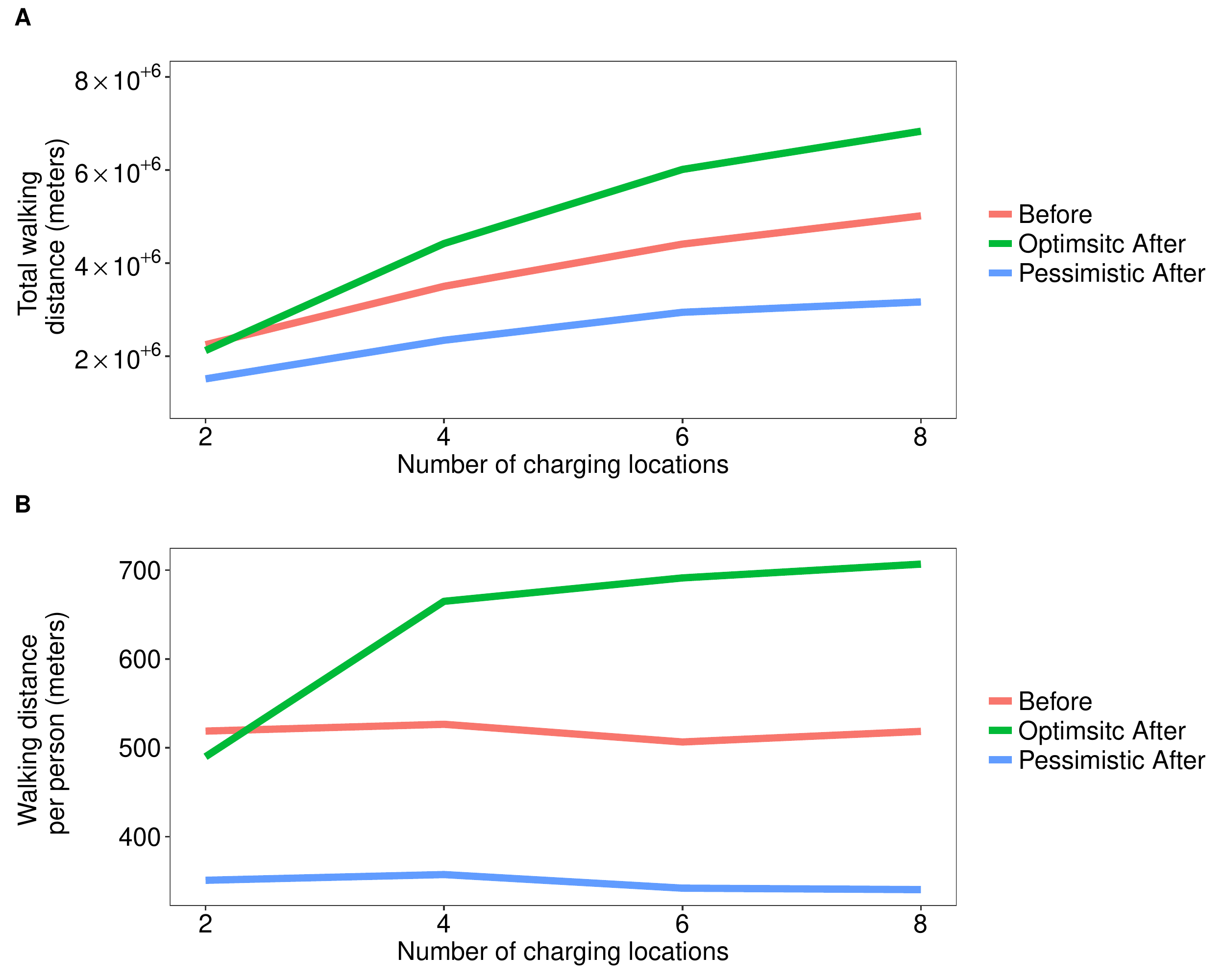}
	\caption {A) Total walking distance and B) walking distance per capita for people that have access to public EV charging service when (BEV,PHEV) market shares are (2\%,3\%).}
	\label{Walk23}
\end{figure}

Figs. \ref{Walk12} and \ref{Walk23} compare total walking distance and walking distance per capita among people that have access to public EV charging service in the community before and after installing charging stations. As mentioned earlier, two cases are evaluated, an optimistic case where we assume that people will always choose the farthest available parking lot and a pessimistic case where people will always choose the nearest available parking lot for EV charging. These plots show that increasing number of charging stations in the community can raise total walking distance and walking distance per capita among people that have access to public EV charging stations up to 40\% in an optimistic case. However, the rate of increase in total walking distance and walking distance per capita decreases as more charging stations are installed in the community. This happens as people get closer to the charging stations and their need to walk is reduced. 

\begin{figure}[!t]
	\centering
	\includegraphics[scale = 0.3]{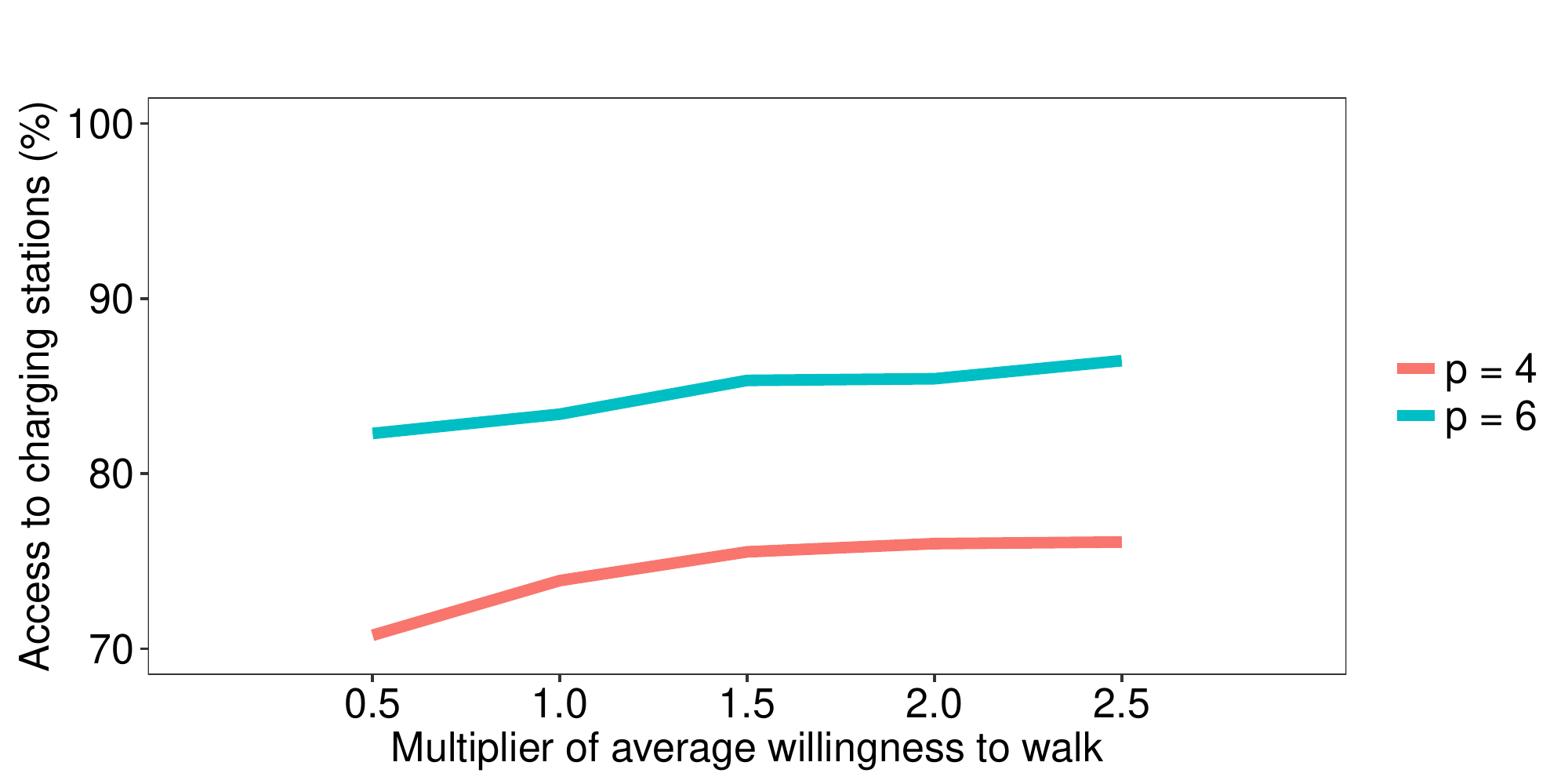}
	\caption {Accessibility for different average of willingness to walk distribution when (BEV,PHEV) market shares are (1\%,2\%).}
	\label{Discussion}
\end{figure}

Another interesting aspect is the relationship between willingness to walk pattern and access to charging stations as young and old communities are expected to have a different level of willingness to walk. Young people tend to walk more while elderly people are not willing to walk long distances. Fig. \ref{Discussion} shows that if the average walking distance preference drops to half, accessibility to public EV charging stations will reduce by 4.23\% and 1.32\% when $p$ = 4 and $p$ = 6, respectively. However, if the average of willingness to walk distribution is doubled, accessibility increases by 2.86\% and 2.43\% when $p$ = 4 and $p$ = 6, respectively. This provides an additional perspective for policy makers, and also indicates the robustness of the model toward any change in willingness to walk pattern in a community.    

\subsection{Value of stochastic solution}

\begin{figure}[!t] 
	\centering
	\includegraphics[scale = 0.3]{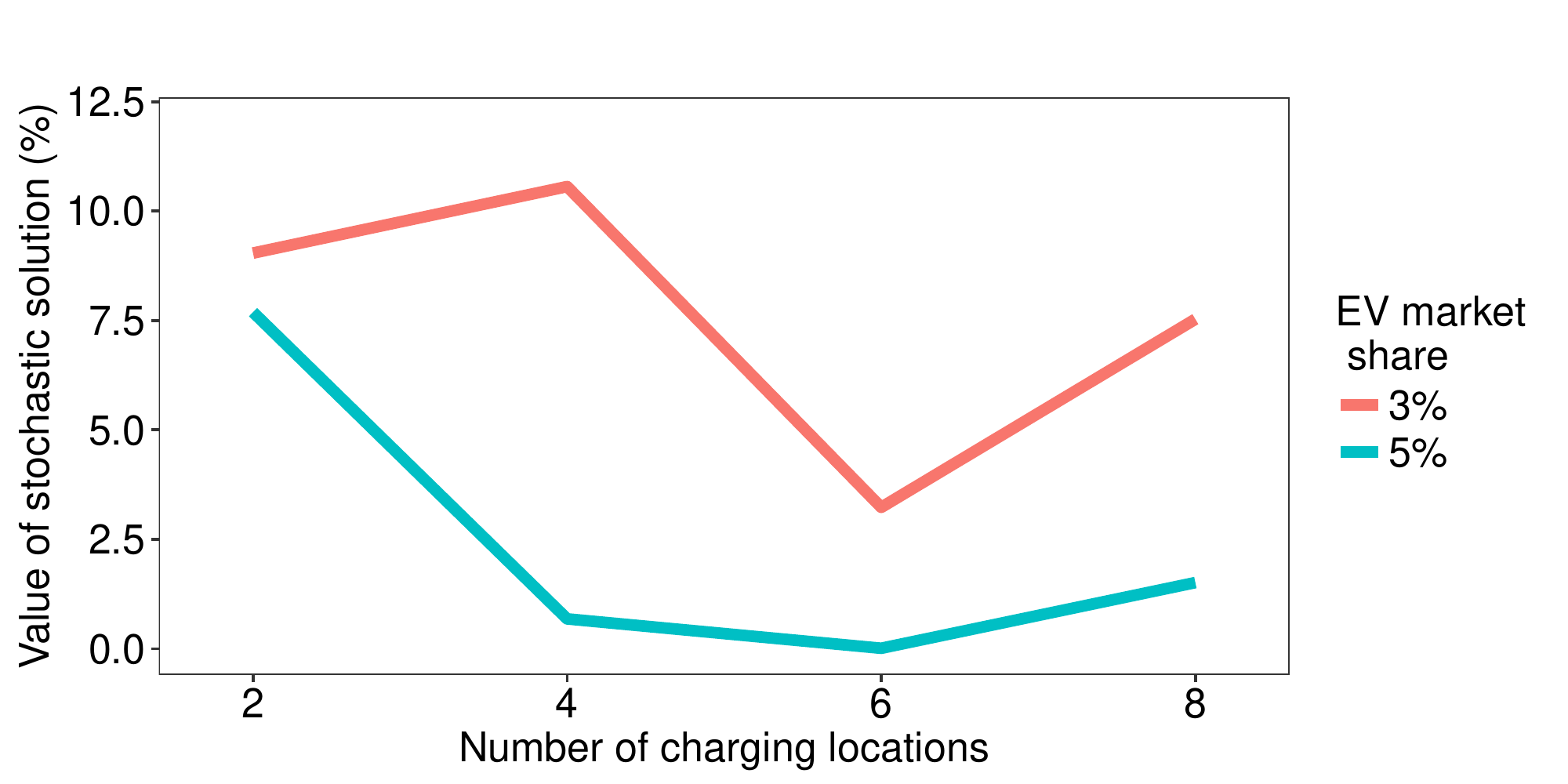}
	\caption {Median of value of stochastic solution for five different runs and different values of $p$ and EV market share.}
	\label{VSS}
\end{figure}

Value of stochastic solution was first introduced by \cite{birge1982value}, and is a standard means to quantify the usefulness of stochastic programming approach. Let the objective value of recourse problem be given as $RP = E_{\Omega}[\varphi(x,z,\omega)]$, and the \textit{expected value problem} is obtained by replacing all random variables in scenarios with their expected values, $EV = \varphi(x,z,\bar{\omega})$, where $\bar{\omega}$ for the demand parameter will be $\sum_{\omega \in \Omega} p_{\omega} d_{\gamma (t),b,s}(\omega)$, $p_{\omega}$ representing probability of occurrence for a scenario $\omega$, and $\sum_{\omega \in \Omega}p_{\omega}=1$. Let $\bar{x},\bar{z}$ represent solutions for $EV$ problem, then the expected result of using expected value solution $(\bar{x},\bar{z})$, is given as $EEV = E_{\Omega}[\varphi(\bar{x},\bar{z},\omega)]$. Then, value of stochastic solution can be defined as $VSS = RP - EEV$. For obtaining VSS, we used the same number of scenarios as in SAA results. Based on five different runs, Fig. \ref{VSS} shows that using stochastic programming brings up to 10.56\% and 7.69\% improvements in accessibility to public EV charging network when EV market share is 3\% and 5\%, respectively.

\section{Conclusion}\label{con}

In this paper, we have presented a two-stage stochastic programming model for public EV charging station network design problem in a community. We considered several uncertainties such as total EV flows, arrival time, dwell time, battery SOC at the time of arrival, charging preference of EV drivers and willingness to walk patterns in estimating demand for public EV charging service. We used sample average approximation method, and for better computational performance, we proposed an effective heuristic that can solve large-scale problems and produce near optimal solutions. On a post analysis, our model presented a number of insights about the design of public EV charging network in an urban/community area. The results show that increasing number of charging stations in the community will improve accessibility to charging service for EV owners but will reduce the utilization level of these stations. Although all charging stations have similar demand patterns but increasing number of charging stations will increase the difference among stations in terms of utilization. While having more charging stations in the community can potentially increase total walking distance and walking distance per capita but the rate of increase in these measures decreases as we install more charging stations. Our model also shows robustness toward any change in willingness to walk pattern of community in the future. We suppose these analogies will provide better insights for a policy maker. Though we have used expected value function for the two-stage model, it will be interesting to see the use of risk-measures for these strategic decisions in the future.

\section{Acknowledgment}

This research was funded by the U.S. DoT through the Transportation Research Center for Livable Communities at Western Michigan University. 

\bibliographystyle{ieeetran}
\bibliography{Bio}

\begin{thebibliography}{10}
\providecommand{\url}[1]{#1}
\csname url@samestyle\endcsname
\providecommand{\newblock}{\relax}
\providecommand{\bibinfo}[2]{#2}
\providecommand{\BIBentrySTDinterwordspacing}{\spaceskip=0pt\relax}
\providecommand{\BIBentryALTinterwordstretchfactor}{4}
\providecommand{\BIBentryALTinterwordspacing}{\spaceskip=\fontdimen2\font plus
\BIBentryALTinterwordstretchfactor\fontdimen3\font minus
  \fontdimen4\font\relax}
\providecommand{\BIBforeignlanguage}[2]{{%
\expandafter\ifx\csname l@#1\endcsname\relax
\typeout{** WARNING: IEEEtran.bst: No hyphenation pattern has been}%
\typeout{** loaded for the language `#1'. Using the pattern for}%
\typeout{** the default language instead.}%
\else
\language=\csname l@#1\endcsname
\fi
#2}}
\providecommand{\BIBdecl}{\relax}
\BIBdecl

\bibitem{berger2015comparison}
D.~J. Berger and A.~D. Jorgensen, ``A comparison of carbon dioxide emissions
  from electric vehicles to emissions from internal combustion vehicles,''
  \emph{Journal of Chemical Education}, vol.~92, no.~7, pp. 1204--1208, 2015.

\bibitem{sierzchula2014influence}
W.~Sierzchula, S.~Bakker, K.~Maat, and B.~van Wee, ``The influence of financial
  incentives and other socio-economic factors on electric vehicle adoption,''
  \emph{Energy Policy}, vol.~68, pp. 183--194, 2014.

\bibitem{upchurch2009model}
C.~Upchurch, M.~Kuby, and S.~Lim, ``A model for location of capacitated
  alternative-fuel stations,'' \emph{Geographical Analysis}, vol.~41, no.~1,
  pp. 85--106, 2009.

\bibitem{frade2011optimization}
I.~Frade, A.~Ribeiro, G.~Goncalves, and A.~Antunes, ``An optimization model for
  locating electric vehicle charging stations in central urban areas,''
  \emph{Journal of Transport Research Record: Journal of the Transportation
  Research Board, DOI}, vol.~10, pp. 2252--12, 2011.

\bibitem{capar2012efficient}
I.~Capar and M.~Kuby, ``An efficient formulation of the flow refueling location
  model for alternative-fuel stations,'' \emph{IIE Transactions}, vol.~44,
  no.~8, pp. 622--636, 2012.

\bibitem{capar2013arc}
I.~Capar, M.~Kuby, V.~J. Leon, and Y.-J. Tsai, ``An arc cover--path-cover
  formulation and strategic analysis of alternative-fuel station locations,''
  \emph{European Journal of Operational Research}, vol. 227, no.~1, pp.
  142--151, 2013.

\bibitem{wang2013locating}
Y.-W. Wang and C.-C. Lin, ``Locating multiple types of recharging stations for
  battery-powered electric vehicle transport,'' \emph{Transportation Research
  Part E: Logistics and Transportation Review}, vol.~58, pp. 76--87, 2013.

\bibitem{baouche2014efficient}
F.~Baouche, R.~Billot, R.~Trigui, and N.-E. El~Faouzi, ``Efficient allocation
  of electric vehicles charging stations: optimization model and application to
  a dense urban network,'' \emph{IEEE Intelligent Transportation Systems
  Magazine}, vol.~6, no.~3, pp. 33--43, 2014.

\bibitem{cavadas2015mip}
J.~Cavadas, G.~H. de~Almeida~Correia, and J.~Gouveia, ``A mip model for
  locating slow-charging stations for electric vehicles in urban areas
  accounting for driver tours,'' \emph{Transportation Research Part E:
  Logistics and Transportation Review}, vol.~75, pp. 188--201, 2015.

\bibitem{he2015deploying}
F.~He, Y.~Yin, and J.~Zhou, ``Deploying public charging stations for electric
  vehicles on urban road networks,'' \emph{Transportation Research Part C:
  Emerging Technologies}, vol.~60, pp. 227--240, 2015.

\bibitem{huang2015optimization}
Y.~Huang and Y.~Zhou, ``An optimization framework for workplace charging
  strategies,'' \emph{Transportation Research Part C: Emerging Technologies},
  vol.~52, pp. 144--155, 2015.

\bibitem{shahraki2015optimal}
N.~Shahraki, H.~Cai, M.~Turkay, and M.~Xu, ``Optimal locations of electric
  public charging stations using real world vehicle travel patterns,''
  \emph{Transportation Research Part D: Transport and Environment}, vol.~41,
  pp. 165--176, 2015.

\bibitem{pan2010locating}
F.~Pan, R.~Bent, A.~Berscheid, and D.~Izraelevitz, ``Locating phev exchange
  stations in v2g,'' in \emph{Smart Grid Communications (SmartGridComm), 2010
  First IEEE International Conference on}.\hskip 1em plus 0.5em minus
  0.4em\relax IEEE, 2010, pp. 173--178.

\bibitem{tan2014stochastic}
J.~Tan and W.-H. Lin, ``A stochastic flow capturing location and allocation
  model for siting electric vehicle charging stations,'' in \emph{17th
  International IEEE Conference on Intelligent Transportation Systems
  (ITSC)}.\hskip 1em plus 0.5em minus 0.4em\relax IEEE, 2014, pp. 2811--2816.

\bibitem{hosseini2015refueling}
M.~Hosseini and S.~MirHassani, ``Refueling-station location problem under
  uncertainty,'' \emph{Transportation Research Part E: Logistics and
  Transportation Review}, vol.~84, pp. 101--116, 2015.

\bibitem{wustochastic}
F.~Wu and R.~Sioshansi, ``A stochastic flow-catching model to optimize the
  location of fast-charging stations with uncertain electric vehicle flows.''

\bibitem{bergantino2014drives}
A.~S. Bergantino, C.~Capozza, A.~De~Carlo, and A.~Morone, ``What drives parking
  choices: a laboratory experiment,'' 2014.

\bibitem{brooker2015identification}
R.~P. Brooker and N.~Qin, ``Identification of potential locations of electric
  vehicle supply equipment,'' \emph{Journal of Power Sources}, vol. 299, pp.
  76--84, 2015.

\bibitem{krumm2012people}
J.~Krumm, ``How people use their vehicles: Statistics from the 2009 national
  household travel survey,'' SAE Technical Paper, Tech. Rep., 2012.

\bibitem{zhong2008studying}
M.~Zhong, J.~Hunt, and X.~Lu, ``Studying differences of household weekday and
  weekend activities: a duration perspective,'' \emph{Transportation Research
  Record: Journal of the Transportation Research Board}, no. 2054, pp. 28--36,
  2008.

\bibitem{ozdemir2015distributed}
A.~Ozdemir \emph{et~al.}, ``Distributed storage capacity modelling of ev
  parking lots,'' in \emph{2015 9th International Conference on Electrical and
  Electronics Engineering (ELECO)}.\hskip 1em plus 0.5em minus 0.4em\relax
  IEEE, 2015, pp. 359--363.

\bibitem{carley2013intent}
S.~Carley, R.~M. Krause, B.~W. Lane, and J.~D. Graham, ``Intent to purchase a
  plug-in electric vehicle: A survey of early impressions in large us cites,''
  \emph{Transportation Research Part D: Transport and Environment}, vol.~18,
  pp. 39--45, 2013.

\bibitem{EPA2016midterm}
``Midterm evaluation of light-duty vehicle greenhouse gas emission standards
  and corporate average fuel economy standards for model years 2022-2025,''
  \emph{Environmental Protection Agency}, 2016.

\bibitem{yang2012walking}
Y.~Yang and A.~V. Diez-Roux, ``Walking distance by trip purpose and population
  subgroups,'' \emph{American journal of preventive medicine}, vol.~43, no.~1,
  pp. 11--19, 2012.

\bibitem{mak1999monte}
W.-K. Mak, D.~P. Morton, and R.~K. Wood, ``Monte carlo bounding techniques for
  determining solution quality in stochastic programs,'' \emph{Operations
  Research Letters}, vol.~24, no.~1, pp. 47--56, 1999.

\bibitem{tsiligirides1984heuristic}
T.~Tsiligirides, ``Heuristic methods applied to orienteering,'' \emph{Journal
  of the Operational Research Society}, vol.~35, no.~9, pp. 797--809, 1984.

\bibitem{vergis2015comparison}
S.~Vergis and B.~Chen, ``Comparison of plug-in electric vehicle adoption in the
  united states: A state by state approach,'' \emph{Research in Transportation
  Economics}, vol.~52, pp. 56--64, 2015.

\bibitem{dai2011optimization}
M.~Dai, J.~Zheng, M.~Zhang, and W.~Wang, ``Optimization of electric vehicle
  charging capacity in a parking lot for reducing peak and filling valley in
  power grid,'' in \emph{Advanced Power System Automation and Protection
  (APAP), 2011 International Conference on}, vol.~2.\hskip 1em plus 0.5em minus
  0.4em\relax IEEE, 2011, pp. 1501--1506.

\bibitem{luo2013optimal}
Z.~Luo, Z.~Hu, Y.~Song, Z.~Xu, and H.~Lu, ``Optimal coordination of plug-in
  electric vehicles in power grids with cost-benefit analysis—part ii: A case
  study in china,'' \emph{IEEE Transactions on Power Systems}, vol.~28, no.~4,
  pp. 3556--3565, 2013.

\bibitem{lutsey2015assessment}
N.~Lutsey, S.~Searle, S.~Chambliss, and A.~Bandivadekar, ``Assessment of
  leading electric vehicle promotion activities in united states cities,''
  \emph{International Council on Clean Transportation}, 2015.

\bibitem{chang2012financial}
D.~Chang, D.~Erstad, E.~Lin, A.~F. Rice, C.~T. Goh, A.~Tsao, and J.~Snyder,
  ``Financial viability of non-residential electric vehicle charging
  stations,'' \emph{Luskin Center for Innovation: Los Angeles, CA, USA}, 2012.

\bibitem{marshall2013evaluation}
W.~E. Marshall, ``An evaluation of livability in creating transit-enriched
  communities for improved regional benefits,'' \emph{Research in
  Transportation Business \& Management}, vol.~7, pp. 54--68, 2013.

\bibitem{litman2004economic}
T.~Litman, ``Economic value of walkability,'' \emph{World Transport Policy \&
  Practice}, vol.~10, no.~1, pp. 5--14, 2004.

\bibitem{birge1982value}
J.~R. Birge, ``The value of the stochastic solution in stochastic linear
  programs with fixed recourse,'' \emph{Mathematical programming}, vol.~24,
  no.~1, pp. 314--325, 1982.

\end{thebibliography}
\newcommand{\BIBdecl}{\setlength{\itemsep}{0.1 em}}
\begin{IEEEbiography}[{\includegraphics[width=1in,height=1.25in]{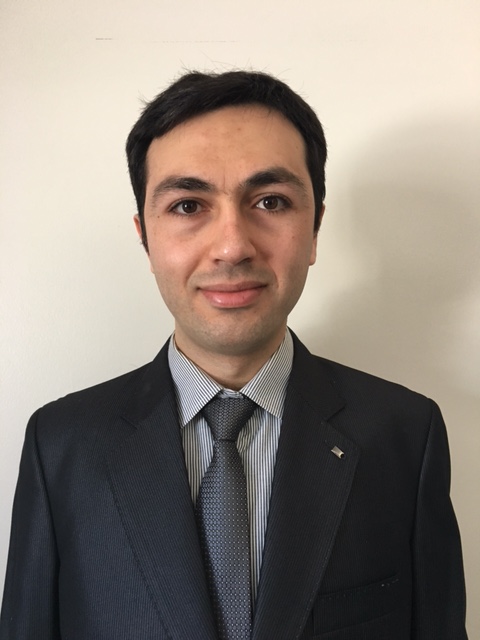}}]{Sina Faridimehr}
	received his B.S. and M.S. degrees in industrial engineering from Sharif University of Technology, Tehran, Iran, in 2009 and 2011 respectively. He is currently working toward his Ph.D. degree in industrial and systems engineering with Wayne State University, Detroit, Michigan. His research interests are in service operations management, and data-driven decision making.
\end{IEEEbiography}

\begin{IEEEbiography}[{\includegraphics[width=1in,height=1.25in]{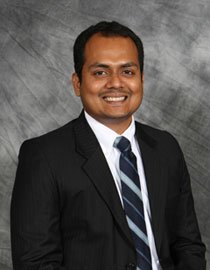}}]{Saravanan Venkatachalam}
	received his M.S. and Ph.D. in industrial and systems engineering from Texas A\&M University, College Station. in 2003 and 2014 respectively. He is currently an Assistant Professor of Industrial and Systems Engineering at Wayne State University. His research interests are in stochastic programming, large scale optimization, and discrete event modeling and simulation. Applications of interest include supply chain management, healthcare, pricing and revenue management, and energy management.
\end{IEEEbiography}

\begin{IEEEbiography}[{\includegraphics[width=1in,height=1.25in]{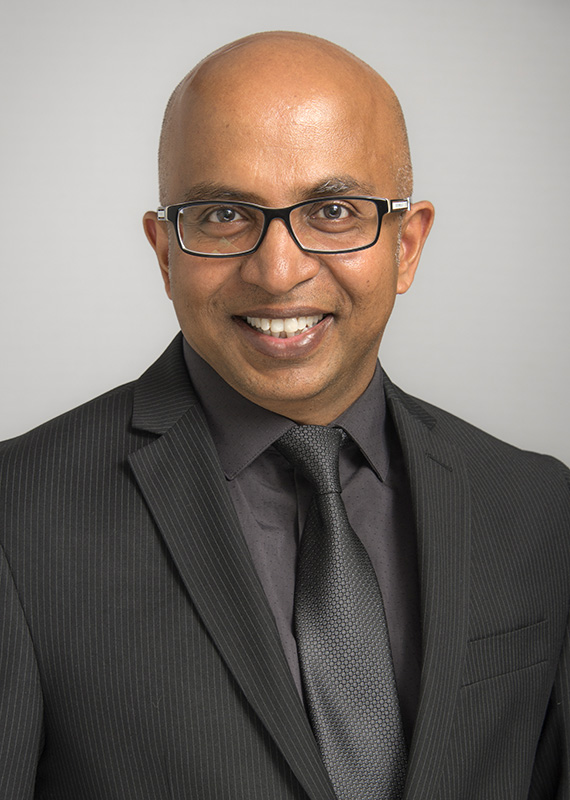}}]{Ratna Babu Chinnam}
	received his  M.S. and Ph.D. degrees in industrial engineering from Texas Tech University (U.S.A.) in 1990 and 1994, respectively. He is currently an Professor of Industrial and Systems Engineering and director of Big Data and Business Analytics Group at Wayne State University. His research interests include  business analytics, big data, supply chain management, freight logistics, operations management, sustainability, healthcare systems engineering, and smart engineering systems.
\end{IEEEbiography}

\end{document}